\documentclass{article}
\usepackage{amsmath,array,amssymb,latexsym, theorem}

\newtheorem{theorem}{Theorem}[section]
\newtheorem{lemma}[theorem]{Lemma}
\newtheorem{proposition}[theorem]{Proposition}

\newtheorem{definition}[theorem]{Definition}
\newtheorem{example}[theorem]{Example}

\newtheorem{remark}[theorem]{Remark}

\numberwithin{equation}{section}

\begin{document}
\noindent
\vspace*{12mm}
\begin{center}
\large{\textbf{\hspace*{5mm}The Geometry of Toric Hyperk\"ahler Varieties} \vspace{12mm} \newline Hiroshi Konno\footnote{Supported in part by the Grant-in-Aid for Scientific Research (C), No. 15540062. \\ \hspace*{4.5mm} 2000 Mathematics Subject Classification: Primary 53C26; Secondary 53D20, 14L24} }
\end{center}
\vspace{6mm}
\begin{abstract}
In this survey article we describe the geometry of toric hyperk\"ahler varieties, which are hyperk\"ahler quotients of the quaternionic vector spaces by tori. 
In particular, we discuss the Betti numbers, the cohomology ring, and variation of hyperk\"ahler structures  of these spaces with many improved results and proofs. 
\end{abstract}
\vspace{7mm}
\section{Introduction}

Let a compact Lie group $G$ act on a hyperk\"ahler manifold $M$, preserving its hyperk\"ahler structure. 
Then we expect to define a quotient space of $M$ by the `quaternionification' of $G$ so that the quotient space 
is also hyperk\"ahler. Although we do not have a notion of `quaternionification' of Lie groups, 
a hyperk\"ahler quotient was introduced as such a kind of things.  
In \cite{HKLR} a hyperk\"ahler quotient was defined as an analogue of a quotient space in symplectic geometry.  
The geometry of hyperk\"ahler quotients is a mixture of geometry of hyperk\"ahler manifolds and symplectic quotients. 

A hyperk\"ahler structure is one of basic geometric structures in Riemannian geometry. 
A hyperk\"ahler manifold is defined to be a Riemannian manifold whose holonomy group is contained in $Sp(\frac{n}{4})$, which is a group on Berger's list of Riemannian holonomy groups \cite{Ber}.  
In other words, it is a Riemannian manifold with three complex structures, 
which satisfies the quaternion relation, and the Riemannian metric is K\"ahler with respect to each complex structure. 
Thus a hyperk\"ahler structure is rather special. So hyperk\"ahler manifolds have quite rich properties. 

On the other hand, symplectic quotients also have rich properties. 
There are many techniques to investigate the topology of compact symplectic quotients. 
In particular, Kirwan proved the surjectivity of the Kirwan map, which is a map from the equivariant cohomology of the original space to the ordinary cohomology of its symplectic quotient. 
This enables one to compute the Betti numbers of symplectic quotients \cite{Ki}. 
Moreover, in many cases a symplectic quotient can be identified with a quotient in geometric invariant theory (GIT for short). 
Due to this identification, one can study these quotients from viewpoints of both symplectic and algebraic geometry. 
Most fundamental examples are semi-projective toric varieties, which are symplectic quotients of the complex vector spaces by tori. 

In this article we describe the geometry of toric hyperk\"ahler varieties, introduced by Bielawski and Dancer \cite{BD}. 
A toric hyperk\"ahler variety is defined to be a hyperk\"ahler quotient of the quaternionic vector space by a torus. 
It is a non-compact hyperk\"ahler variety on which a torus acts preserving its hyperk\"ahler structure. 
It is a hyperk\"ahler analogue of, but different from an ordinary toric variety. 
We explain fundamental results of toric hyperk\"ahler varieties. 
We give many improved statements and proofs, which cannot be found in literatures. 
We refer the reader to \cite{P} for further results for toric hyperk\"ahler varieties. 

In Sections 2 we define toric hyperk\"ahler varieties and introduce basic notations. 

In Section 3 we prove basic properties of toric hyperk\"ahler varieties. 
First we identify hyperk\"ahler quotients with GIT quotients in our setting. 
Then we give characterizations of semi-stable points and closed orbits, 
which are very useful in this article. 
A toric hyperk\"ahler variety has a parameter, 
which is shown to correspond to the period of its hyperk\"ahler structure. 
Then we introduce natural morphisms between toric hyperk\"ahler varieties with different parameters, 
influenced by \cite{Kr,N1}.  
Then we show that a quotient for a generic parameter is an orbifold, 
and that its topology is independent of the parameter. 
A quotient for a certain parameter may have singularities which are worse than orbifold singularities.

In Section 4 we study the topology of toric hyperk\"ahler varieties for generic parameters. 
Although a projective toric variety is determined by a convex polytope, 
a toric hyperk\"ahler variety is determined by an arrangement of hyperplanes. 
We describe the Betti numbers  and the cohomology ring  in terms of the arrangement of hyperplanes, 
based on the works \cite{BD} and \cite{Ko1} respectively. 
We also describe another computation due to \cite{HS}, where a toric hyperk\"ahler variety is shown to be 
homotopy equivalent to a certain non-compact ordinary toric variety. 
We give a new proof of this homotopy equivalence. 
We also give another useful presentation of the cohomology ring, which is a simplification of the one in \cite{Ko1}.
In the case of hyperk\"ahler quotients we have an analogue of the Kirwan map, 
but it is not yet known whether it is surjective or not in general. The above computation shows 
that the analogue of the Kirwan map is surjective in the case of toric hyperk\"ahler varieties. 

In Section 5 we give two examples, which illustrate the results in other sections.  

In Section 6 we describe how a toric hyperk\"ahler variety changes according to its period. 
In particular, we show that the so called wall-crossing phenomena are described as Mukai flops, 
which are special bimeromorphic transformations for hyperk\"ahler varieties. 
We give a proof of this result using the technical tools developed in Section 3. 
This simplifies the proofs in \cite{Ko3}. 

The author would like to thank the organizers for the exciting workshop 
`Toric Topology Conference 2006' at Osaka City University.
\section{Toric hyperk\"ahler varieties}
Let $\mathbb{H}^N$ be the quaternionic vector space with the 
standard hyperk\"ahler structure $(g, I_1,I_2,I_3)$. Here $g$ is 
the standard Riemannian metric on $\mathbb{H}^N \cong \mathbb{R}^{4N}$, 
and the three complex structures $I_1,I_2,I_3$ are defined by 
the left multiplication by $I_1,I_2,I_3$, which are the standard 
generators of the quaternion $\mathbb{H}$. From now on we identify 
$\xi \in \mathbb{H}^N$ with $(z,w) \in \mathbb{C}^N \times \mathbb{C}^N$ 
by $\xi =z + w I_2$. The real torus
$T^N=\{ \zeta=(\zeta_1, \dots, \zeta_N) \in 
\mathbb{C}^N ~|~ \vert \zeta_i \vert =1 \}$ acts on 
$\mathbb{H}^N$ by $(z,w)\zeta=(z\zeta,w\zeta^{-1})$, 
preserving its hyperk\"ahler structure. 
Let $K$ be a connected subtorus of $T^N$ with the Lie algebra $k \subset t^N$. 
Then we have the following exact sequences;
\[
\begin{array}{ccccccccc}
0 & \longrightarrow & k & \overset{\iota}{\longrightarrow} & t^N 
& \overset{\pi}{\longrightarrow} & t^n & \longrightarrow & 0, \\
0 & \longleftarrow & k^* & \overset{\iota^*}{\longleftarrow} & (t^N)^* 
& \overset{\pi^*}{\longleftarrow} & (t^n)^* & \longleftarrow & 0, 
\end{array}
\]
where $t^n =t^N /k$ is the Lie algebra of the quotient torus $T^n =T^N /K$. 
The action of $K$ on $\mathbb{H}^N$ admits a hyperk\"ahler moment map
\[
\mu=(\mu_{1},\mu_{2},\mu_{3}) \colon 
\mathbb{H}^N \to k^* \otimes \mathbb{R}^3, 
\]
which is given by
\begin{align*}
& \mu_{1}(z,w) = \pi \sum_{i=1}^N (\vert z_i \vert^2 - \vert w_i \vert^2) 
  \iota^* u_i, \\
& (\mu_{2}+\sqrt{-1}\mu_{3})(z,w) 
 = -2 \pi \sqrt{-1} \sum_{i=1}^N z_i w_i \iota^* u_i, 
\end{align*}
where $\{u_1, \dots, u_N \} \subset (t^N)^*$ is the dual basis 
of the standard basis $\{X_1, \dots, X_N \} \subset t^N$.
Note that the complex moment map $\mu_{\mathbb{C}}= \mu_{2}+\sqrt{-1}\mu_{3} 
\colon \mathbb{H}^N \to k_{\mathbb{C}}^*$ 
is holomorphic with respect to $I_1$, where $k_{\mathbb{C}}$ is the 
complexification of the Lie algebra $k$. Then Bielawski and Dancer 
introduced the following. 

\begin{definition}\label{toric} 
\rm{\cite{BD} A \textit{toric hyperk\"ahler variety} $X(\alpha,\beta)$ is a hyperk\"ahler quotient 
$\mu^{-1}(\alpha,\beta) /K$ for $(\alpha,\beta)  \in k^* \times k_{\mathbb{C}}^*$.}
\end{definition}
\begin{remark}
\rm{Toric hyperk\"ahler varieties are sometimes called `hypertoric varieties' \cite{P}. 
They are just the same objects. }
\end{remark}
A smooth part of $X(\alpha,\beta)$ is a $4n$-dimensional 
hyperk\"ahler manifold. We denote its hyperk\"ahler structure by 
$(g, I_{1}, I_{2}, I_{3})$ again. 
The torus $T^n =T^N /K$ acts on $X(\alpha,\beta)$, preserving 
its hyperk\"ahler structure. Note that toric hyperk\"ahler varieties 
are not toric varieties in the ordinary sense, but the corresponding objects 
in hyperk\"ahler geometry to them. 

Bielawski gave an intrinsic characterization of toric hyperk\"ahler manifolds as follows. 

\begin{theorem} {\rm \cite{Bi}}
Let $M$ be a complete connected hyperk\"ahler manifold of dimension 
$4n$. Suppose that the torus $T^n$ acts on $M$ effectively, preserving 
the hyperk\"ahler structure, and that $M$ has Euclidean volume growth. 
Then $M$ is isomorphic to one of toric hyperk\"ahler manifolds in the sense 
of Definition \ref{toric} as a hyperk\"ahler $T^n$-manifold.
\end{theorem}

\section{Holomorphic descriptions}
A toric hyperk\"ahler variety $X(\alpha,\beta)=\mu_{1}^{-1}(\alpha) \cap \mu_{\mathbb{C}}^{-1}(\beta)/K$ can be considered as a symplectic quotient of $\mu_{\mathbb{C}}^{-1}(\beta)$ by the torus $K$. 
In many cases symplectic quotients can be identified with GIT quotients. 
In this section we describe this identification and apply it to prove basic properties of toric hyperk\"ahler varieties. 

\subsection{Symplectic and GIT quotients}
Let $t^N_{\mathbb{Z}}=\sum_{i=1}^N \mathbb{Z}X_i$ and $(t^N_{\mathbb{Z}})^*=\sum_{i=1}^N \mathbb{Z}u_i$ be the lattices in $t^N$ and $(t^N)^*$ respectively. 
We also set $k_{\mathbb{Z}}=k \cap t^N_{\mathbb{Z}}$ and $k_{\mathbb{Z}}^*= Hom_{\mathbb{Z}}(k_{\mathbb{Z}}, \mathbb{Z})$. 
The element $\alpha \in k_{\mathbb{Z}}^*$ induces the character $\chi_{\alpha} \colon K_{\mathbb{C}} \to \mathbb{C}^{\times}$, where $K_{\mathbb{C}}$ is the complexfication of $K$. 
Let $L^{\otimes m}=\mathbb{H}^N \times \mathbb{C}$ be the trivial holomorphic line bundle on which $K_{\mathbb{C}}$ acts by 
$$
((z,w), v)_m \zeta= ((z \zeta, w \zeta^{-1}), v \chi_{\alpha}(\zeta)^m)_m.
$$

Let us consider the GIT quotient of the affine variety $\mu_{\mathbb{C}}^{-1}(\beta)$ by $K_{\mathbb{C}}$ with respect to the linearization induced by $\alpha \in k_{\mathbb{Z}}^*$.
Recall that a point $(z,w) \in \mu_{\mathbb{C}}^{-1}(\beta)$ is $\alpha$-{\it semi-stable} if and only if there exists $m \in \mathbb{Z}_{>0}$ and a polynomial $f(p, q)$, where $p,q \in \mathbb{C}^N$, such that $f((p,q)\zeta)=f(p,q)\chi_{\alpha}(\zeta)^m$ for any $\zeta \in K_{\mathbb{C}}$ and $f(z,w) \ne 0$. 
We denote the set of $\alpha$-semi-stable points in $\mu_{\mathbb{C}}^{-1}(\beta)$ by $\mu_{\mathbb{C}}^{-1}(\beta)^{\alpha-ss}$. 
Then we have a categorical quotient $\phi \colon \mu_{\mathbb{C}}^{-1}(\beta)^{\alpha-ss} \to 
\mu_{\mathbb{C}}^{-1}(\beta)^{\alpha-ss}// K_{\mathbb{C}}$. 
Here we summarize the properties of categorical quotients. We refer the reader to \cite{Do,MFK} for the proof. 

\begin{lemma}\label{cat}
For any point $p \in \mu_{\mathbb{C}}^{-1}(\beta)^{\alpha-ss}// K_{\mathbb{C}}$, the fiber $\phi^{-1}(p)$ consists of finitely many $K_{\mathbb{C}}$-orbits. 
Moreover, each fiber contains the unique closed $K_{\mathbb{C}}$-orbits in $\mu_{\mathbb{C}}^{-1}(\beta)^{\alpha-ss}$. 
Thus the categorical quotient $\mu_{\mathbb{C}}^{-1}(\beta)^{\alpha-ss}// K_{\mathbb{C}}$ can be identified with a set of closed $K_{\mathbb{C}}$-orbits in $\mu_{\mathbb{C}}^{-1}(\beta)^{\alpha-ss}$.  
\end{lemma}

Define a Hermitian metric on the holomorphic line bundle $L^{\otimes m}=\mathbb{H}^N \times \mathbb{C}$ by $||((z,w), v)_m||=|v| e^{-\frac{\pi}{2}m(|z|^2 +|w|^2)}$. 
It induces the Chern connection $\nabla^{\otimes m}$ on $L^{\otimes m}$. 
Note that the first Chern form $c_1(\nabla^{\otimes m})=m \omega_1$, where $\omega_1 \in \Omega^2(\mathbb{H}^N)$ is the standard K\"ahler form with respect to $I_1$. 
Since the action of the real torus $K$ on $L^{\otimes m}$ preserves the holomorphic structure and the Hermitian metric, it also preserves the Chern connection $\nabla^{\otimes m}$. 
Let us introduce the following function. 
Fix $(z,w) \in \mathbb{H}^N$. 
Then define $l_{(z,w)} \colon K_{\mathbb{C}} \to \mathbb{R}$ by 
$$ 
l_{(z,w)}(\zeta)=\frac{1}{2\pi} \log \frac{||((z,w), v)_{-1} \zeta||}{||((z,w), v)_{-1}||}.
$$
In particular, for $X \in k$ and $t \in \mathbb{R}$, we have 
\begin{equation}\label{fun}\begin{split}
l_{(z,w)}& (\mathrm{Exp} \sqrt{-1}tX) =\langle \alpha, tX \rangle  \\
& + \frac{1}{4}\sum_{i=1}^N |z_i |^2 e^{-4 \pi \langle u_i, tX \rangle} + 
\frac{1}{4}\sum_{i=1}^N |w_i |^2 e^{4 \pi \langle u_i, tX \rangle} - \frac{1}{4}(|z|^2+|w|^2), 
\end{split}\end{equation} 
where $\mathrm{Exp} \colon k_{\mathbb{C}} \to K_{\mathbb{C}}$ is the exponential map. Then we have the following. 
\begin{lemma}\label{moment}
$(1)$ $l_{(z,w)}$ is a convex function. 
\newline
$(2)$ $\zeta$ is a critical point if and only if $\mu_1((z,w)\zeta)= \alpha$.
\newline
$(3)$ All critical points are minima of $l_{(z,w)}$.
\newline
$(4)$ The following $(i)$, $(ii)$ and $(iii)$ are equivalent:
\newline
\hspace*{5mm}$(i)$ $(z,w) \in \mu_{\mathbb{C}}^{-1}(\beta)^{\alpha-ss}$ and the $K_{\mathbb{C}}$-orbit through $(z,w)$ is closed in $\mu_{\mathbb{C}}^{-1}(\beta)^{\alpha-ss}$.
\newline
\hspace*{5mm}$(ii)$ $(z,w)K_{\mathbb{C}} \cap \mu_1^{-1}(\alpha) \ne \emptyset$. 
\newline
\hspace*{5mm}$(iii)$ $l_{(z,w)}$ attains a minimum. 
\end{lemma}
\textit{Proof.} A direct computation shows 
\begin{equation}
\begin{split}
& \frac{d}{dt}l_{(z,w)} (\mathrm{Exp} \sqrt{-1}tX)|_{t=t_0}=
\langle \alpha -\mu_1((z,w)\mathrm{Exp} \sqrt{-1}t_0X),X \rangle, \\
& \frac{d^2}{dt^2}l_{(z,w)} (\mathrm{Exp} \sqrt{-1}tX)|_{t=t_0}= 
|X^*|^2_{\text{at $(z,w)\mathrm{Exp} \sqrt{-1}t_0X$} },
\end{split}
\end{equation}
where $X^*$ is a vector field on $\mu_{\mathbb{C}}^{-1}(\beta)$ induced by $X \in k$. 
Then it is easy to see $(1)$, $(2)$ and $(3)$.  
We refer the reader to \cite{N2} for the proof of $(4)$.  $\hfill \Box$
\vspace{3mm} \newline 
By using Lemmas \ref{cat} and \ref{moment}, we can show the following. 
\begin{proposition}\label{GIT}
If we fix $\alpha \in k_{\mathbb{Z}}^*$, then the natural map $\iota \colon X(\alpha, \beta) \to 
\mu_{\mathbb{C}}^{-1}(\beta)^{\alpha-ss}// K_{\mathbb{C}}$ is a homeomorphism. 
\end{proposition}
Thus we can identify the symplectic quotient $(X(\alpha, \beta),I_1)$ with the GIT quotient $\mu_{\mathbb{C}}^{-1}(\beta)^{\alpha-ss}// K_{\mathbb{C}}$. 
This principle was established in \cite{KN,MFK} in the case of projective quotients. 
The above version is proved in \cite{N2}.

\subsection{Characterization of semi-stable points and closed orbits}
The following technical lemma will be very useful in this article.  
\begin{lemma}\label{criterion}
Suppose that $\alpha \in k^*_{\mathbb{Z}}$. 
\newline
$(1)$ A point $(z,w) \in \mu_{\mathbb{C}}^{-1}(\beta)$ is 
$\alpha$-{\it semi-stable} if and only if 
\begin{equation}\label{ss}
\alpha \in \sum_{i=1}^N \mathbb{R}_{\ge 0}|z_i|^2\iota^{*}u_i
+\sum_{i=1}^N \mathbb{R}_{\ge 0}|w_i|^2 (-\iota^{*}u_i). 
\end{equation}
\newline
$(2)$ Suppose $(z,w) \in \mu_{\mathbb{C}}^{-1}(\beta)^{\alpha-ss}$. Then the $K_{\mathbb{C}}$-orbit 
through $(z,w)$ is closed in $\mu_{\mathbb{C}}^{-1}(\beta)^{\alpha-ss}$ if and only if 
\begin{equation}\label{closed}
\alpha \in \sum_{i=1}^N \mathbb{R}_{> 0}|z_i|^2\iota^{*}u_i
+\sum_{i=1}^N \mathbb{R}_{> 0}|w_i|^2 (-\iota^{*}u_i).
\end{equation}
\end{lemma}
\textit{Proof.} 
$(1)$ Suppose $(z,w) \in \mu_{\mathbb{C}}^{-1}(\beta)^{\alpha -ss}$. Then there exists $m \in \mathbb{Z}_{>0}$ 
and a polynomial $f(p_1, \dots, p_N, q_1, \dots q_N)$ such that $f((p,q)\zeta)=f(p,q)\chi_{\alpha}(\zeta)^m$ 
for any $\zeta \in K_{\mathbb{C}}$ and $f(z,w) \ne 0$. So there exists a monomial 
$f_0(p,q)=\prod_{i=1}^Np_i^{a_i}\prod_{i=1}^Nq_i^{b_i}$, where $a_i,b_i \in \mathbb{Z}_{\ge 0}$, such that 
$f_0((p,q)\zeta)=f_0(p,q)\chi_{\alpha}(\zeta)^m$ and $f_0(z,w) \ne 0$. The second condition implies that
$a_i=0$ if $z_i=0$, and $b_i=0$ if $w_i=0$. Moreover, the first condition implies 
$m \alpha=\sum_{i=1}^Na_i \iota^* u_i + \sum_{i=1}^Nb_i (-\iota^* u_i)$. Thus we proved (\ref{ss}). 

Suppose that (\ref{ss}) holds. Then there exists $m \in \mathbb{Z}_{>0}$ and 
$a_i,b_i \in \mathbb{Z}_{\ge 0}$ such that 
$m \alpha=\sum_{i=1}^Na_i \iota^* u_i + \sum_{i=1}^Nb_i (-\iota^* u_i)$ and 
such that $a_i=0$ if $z_i=0$, and $b_i=0$ if $w_i=0$. Then the monomial $f_0(p,q)=\prod_{i=1}^Np_i^{a_i}\prod_{i=1}^Nq_i^{b_i}$ satisfies 
$f_0((p,q)\zeta)=f_0(p,q)\chi_{\alpha}(\zeta)^m$ and $f_0(z,w) \ne 0$. Therefore we have 
$(z,w) \in \mu_{\mathbb{C}}^{-1}(\beta)^{\alpha -ss}$. 
\newline
$(2)$ Suppose that the $K_{\mathbb{C}}$-orbit through $(z,w)$ is closed in $\mu_{\mathbb{C}}^{-1}(\beta)^{\alpha-ss}$. Then, by Lemma \ref{moment}, there exists $\zeta \in K_{\mathbb{C}}$ 
such that $\alpha=\mu_1((z,w)\zeta)$. Moreover, we have 
$\mu_1((z,w)\zeta) \in \sum_{i=1}^N \mathbb{R}_{> 0}|z_i|^2\iota^{*}u_i
+\sum_{i=1}^N \mathbb{R}_{> 0}|w_i|^2 (-\iota^{*}u_i)$. This implies (\ref{closed}).

Suppose that (\ref{closed}) holds. By Lemma \ref{moment}, we have only to show 
that $l_{(z,w)}$ attains a minimum. We claim that 
\begin{equation}\label{claim1}
\lim_{t \to \infty}l_{(z,w)}(\mathrm{Exp} \sqrt{-1}tX) \ne \infty ~~
\text{if and only if}~~X \in k_{(z,w)}, 
\end{equation}
where $k_{(z,w)}$ is the Lie algebra of the isotropy subgroup 
$K_{(z,w)}$ of $K$ at $(z,w)$. In fact, by (\ref{fun}), the condition 
$\displaystyle \lim_{t \to \infty}l_{(z,w)}(\mathrm{Exp} \sqrt{-1}tX) \ne \infty$ implies the following 
conditions $(i)$ and $(ii)$;
\newline
\hspace*{5mm}$(i)$ $\langle \alpha, X \rangle \le 0$,
\newline
\hspace*{5mm}$(ii)$ $\langle u_i, X \rangle \ge 0$ if $z_i \ne 0$, and 
$\langle u_i, X \rangle \le 0$ if $w_i \ne 0$.
\newline
Due to the assumption (\ref{closed}), $(ii)$ implies $\langle \alpha, X \rangle \ge 0$. 
Due to $(i)$, we conclude that $\langle u_i, X \rangle = 0$ if $z_i \ne 0$, and that 
$\langle u_i, X \rangle = 0$ if $w_i \ne 0$. This implies that $X \in k_{(z,w)}$. 
On the other hand, it is easy to see that, if $X \in k_{(z,w)}$, then $l_{(z,w)}(\mathrm{Exp} \sqrt{-1}tX)$ is independent of $t \in \mathbb{R}$. 
Thus we proved the claim (\ref{claim1}). 
Then, using the convexity of $l_{(z,w)}$, we can easily show that 
$l_{(z,w)}$ attains a minimum. $\hfill \Box$
\vspace{3mm} \newline 
By this lemma the notion of $\alpha$-semi-stability is defined for any $\alpha \in k^*$. 
Namely, a point $(z,w) \in \mu_{\mathbb{C}}^{-1}(\beta)$ is defined to be $\alpha$-{\it semi-stable} 
if and only if 
\begin{equation*}
\alpha \in \sum_{i=1}^N \mathbb{R}_{\ge 0}|z_i|^2\iota^{*}u_i
+\sum_{i=1}^N \mathbb{R}_{\ge 0}|w_i|^2 (-\iota^{*}u_i) 
\end{equation*}
for any $\alpha \in k^*$. Moreover, since the above argument works for any $\alpha \in k^*$, 
we can identify a toric hyperk\"ahler variety $(X(\alpha, \beta), I_1)$ with 
$\mu_{\mathbb{C}}^{-1}(\beta)^{\alpha-ss}// K_{\mathbb{C}}$ for any $\alpha \in k^*$. 

Applying this lemma to the case $\alpha=0 \in k^*$, we have 
$\mu_{\mathbb{C}}^{-1}(\beta)^{0-ss}=\mu_{\mathbb{C}}^{-1}(\beta)$. 
So the toric hyperk\"ahler variety $X(0,\beta)= \mu_{\mathbb{C}}^{-1}(\beta)// K_{\mathbb{C}}$ 
is an affine variety.  
\subsection{Natural Morphisms}
The inclusion $\mu_{\mathbb{C}}^{-1}(\beta)^{\alpha-ss} \subset \mu_{\mathbb{C}}^{-1}(\beta)$ 
induces a natural morphism from a GIT quotient $\mu_{\mathbb{C}}^{-1}(\beta)^{\alpha-ss}// K_{\mathbb{C}}$ to 
an affine quotient $\mu_{\mathbb{C}}^{-1}(\beta)// K_{\mathbb{C}}$, which we denote by 
$$\pi_1 \colon (X(\alpha, \beta), I_1) \to (X(0, \beta),I_1).$$ 
If we denote $X(\alpha, \beta)$ by $X(\alpha_1, \alpha_2, \alpha_3)$, where $(\alpha, \beta)=
(\alpha_1, \alpha_2 + \sqrt{-1}\alpha_3)$, then the natural morphism is written 
as $\pi_1 \colon (X(\alpha_1, \alpha_2, \alpha_3), I_1) \to 
(X(0, \alpha_2, \alpha_3),I_1)$.  By rotating the complex structures 
$I_1$, $I_2$ and $I_3$ cyclically, we also have the natural morphisms 
\begin{align*}
& \pi_2 \colon (X(\alpha_1, \alpha_2, \alpha_3), I_2) \to (X(\alpha_1, 0, \alpha_3),I_2), \\ 
& \pi_3 \colon (X(\alpha_1, \alpha_2, \alpha_3), I_3) \to (X(\alpha_1, \alpha_2, 0),I_3).
\end{align*}
The following is a basic property of the natural morphisms, which will be useful in the next section.
\begin{proposition}\label{proper}
The natural morphism $\pi_1 \colon (X(\alpha,\beta),I_1) \to (X(0, \beta),I_1)$ is proper and surjective for any $(\alpha, \beta) \in k^* \times k_{\mathbb{C}}^*$.
\end{proposition}
\textit{Proof.} 
First we prove that $\pi_1$ is proper. Suppose the contrary. 
Namely, there exists a compact set $S \subset X(0,\beta)$ such that $\pi_1^{-1}(S) \subset X(\alpha, \beta)$ is non-compact. 
Then we can take a sequence $\{ p_{(i)} \} \subset \mu^{-1}(\alpha,\beta)$ such that $[p_{(i)}] \in \pi^{-1}(S)$ and $||p_{(i)}|| \to \infty$. 
If we set $p_{(i)}^{\prime}=\frac{p_{(i)}}{||p_{(i)}||} \in \mu^{-1}(\frac{\alpha}{||p_{(i)}||^2},\frac{\beta}{||p_{(i)}||^2})$, 
then we may assume that $\{ p_{(i)}^{\prime} \}$ converges to $p_{(\infty)}^{\prime} \in \mathbb{H}^N$. 
Then we have $|| p_{(\infty)}^{\prime} ||=1$ and $p_{(\infty)}^{\prime} \in \mu^{-1}(0,0)$. 

On the other hand, we can take $q_{(i)} \in \mu^{-1}(0,\beta)$ such that $\pi_1([p_{(i)}])=[q_{(i)}]$. 
Since $\{ [q_{(i)}] \} \subset S$, $\{ q_{(i)} \}$ is bounded in $\mu^{-1}(0,\beta)$. 
If we set $q_{(i)}^{\prime}=\frac{q_{(i)}}{||p_{(i)}||} 
\in \mu^{-1}(0,\frac{\beta}{||p_{(i)}||^2})$, 
then $\{ q_{(i)}^{\prime} \}$ converges to $(0,0) \in \mathbb{H}^N$. 
Let $\pi_1^{\prime} \colon X(\frac{\alpha}{||p_{(i)}||^2},\frac{\beta}{||p_{(i)}||^2}) 
\to X(0,\frac{\beta}{||p_{(i)}||^2})$ be the natural morphism. Since $\pi_1^{\prime}([p_{(i)}^{\prime}])=
[q_{(i)}^{\prime}]$, we have $q_{(i)}^{\prime} \in \overline{p_{(i)}^{\prime}K_{\mathbb{C}}}$. 
Therefore we can take $\{ \zeta_{(i)} \} \subset K_{\mathbb{C}}$ such that $\{ p_{(i)}^{\prime} \zeta_{(i)} \}$ converges to $(0,0)$. 
Then we can show that $\{ p_{(\infty)}^{\prime} \zeta_{(i)} \}$ converges to $(0,0)$. 
So we have $(0,0) \in \overline{p_{(\infty)}^{\prime}K_{\mathbb{C}}}$. 

However, since $p_{(\infty)}^{\prime} \in \mu^{-1}(0,0)$, we have 
$\overline{p_{(\infty)}^{\prime}K_{\mathbb{C}}}=p_{(\infty)}^{\prime}K_{\mathbb{C}}$. 
Therefore we have $p_{(\infty)}^{\prime}=(0,0)$. This contradicts to $|| p_{(\infty)}^{\prime} ||=1$. 

Surjectivity of the map $\pi_1$ will be proved at the end of Subsection \ref{the core}. $\hfill \Box$
\vspace{3mm} \newline 
The above proof is essentially the same as the one in \cite{Kr,N1}, where the similar result is proved 
in the case of quiver varieties. See also \cite{Ao} for the proof.

Let $(k^* \times k_{\mathbb{C}}^*)_{reg}$ be the set of regular values 
of the hyperk\"ahler moment map $\mu=(\mu_{1},\mu_{\mathbb{C}}) \colon 
\mathbb{H}^N \to k^* \otimes \mathbb{R}^3 \cong 
k^* \times k_{\mathbb{C}}^*$. 
Denote the isotropy subgroup of $K$ at $(z,w) \in \mathbb{H}^N$ by 
$K_{(z,w)}$. Set 
\begin{align*}
& \Lambda = \{ K_{(z,w)} ~|~ (z,w) \in \mathbb{H}^N \}, \\
& \Lambda^{(d)} = \{ H \in \Lambda ~|~ \dim H =d \} ~~\text{for}~ 
d=0, \dots, N-n.  
\end{align*}
We write $\Lambda^{(1)}= \{ H^{(1)}_1, \dots, H^{(1)}_l \}$. 
The subspace of codimension one 
$$
W_s = \{ v \in k^* ~|~ \langle v, \text{Lie} H^{(1)}_s\rangle =0 \} \subset k^* 
$$
is called a {\it wall}. 
Note that the wall $W_s$ is spanned by $\{ \iota^* u_i ~|~ \langle \iota^* u_i, \text{Lie} H^{(1)}_s\rangle =0 \}$.
Moreover, we have the following. 
\begin{proposition}\label{indep}
$(1)$ $ (k^* \times k_{\mathbb{C}}^*)_{reg} = 
k^* \times k_{\mathbb{C}}^* \setminus \bigcup_{s=1}^l 
W_s \otimes W_{s \mathbb{C}}$
\newline
$(2)$ If $(\alpha,\beta) \in (k^* \times 
k_{\mathbb{C}}^*)_{reg}$, then $X(\alpha,\beta)$ is an 
orbifold. 
\newline
$(3)$ If $(0,\beta) \in (k^* \times k_{\mathbb{C}}^*)_{reg}$, then  $\pi_1 \colon (X(\alpha,\beta), I_1) \to (X(0,\beta),I_1)$ is a biholomorphic map for any $\alpha \in k^*$. 
\newline
$(4)$ If $(\alpha,\beta) \in (k^* \times k_{\mathbb{C}}^*)_{reg}$, then 
every $K_{\mathbb{C}}$-orbit is closed 
in $\mu_{\mathbb{C}}^{-1}(\beta)^{\alpha -ss}$. So the categorical quotient 
$\mu_{\mathbb{C}}^{-1}(\beta)^{\alpha -ss}//K_{\mathbb{C}}$ is a geometric quotient 
$\mu_{\mathbb{C}}^{-1}(\beta)^{\alpha -ss}/K_{\mathbb{C}}$.
\newline
$(5)$ The topology of $X(\alpha,\beta)$ is independent of $(\alpha,\beta)$ as far as $(\alpha,\beta) \in 
(k^* \times k_{\mathbb{C}}^*)_{reg}$. That is, all the $X(\alpha,\beta)$ are diffeomorphic to each other as far as $(\alpha,\beta) \in (k^* \times k_{\mathbb{C}}^*)_{reg}$. The diffeomorphisms are given by compositions of natural morphisms. 
\end{proposition}
\textit{Proof.} 
$(1)$ Let $f \colon \mathbb{H} \to \mathbb{R} \times \mathbb{C}$ be a map defined by 
$f(a,b)=(\pi (|a|^2 -|b|^2), -2\pi \sqrt{-1}ab)$. We can easily observe that 
$(a,b) \in \mathbb{C} \times \mathbb{C} \cong 
\mathbb{H}$ is a regular point of $f$ if and only if $(a,b) \ne (0,0)$. 
Since $(d \mu_K)_{(z,w)}= \sum_{i=1}^N (df)_{(z_i,w_i)}\otimes \iota^* u_i \in k^* \times k_{\mathbb{C}}^*$, 
$(z,w) \in \mathbb{H}^N$ is a critical point of $\mu_K$ if and only if 
span$\{ \iota^* u_i ~|~ (z_i, w_i) \ne (0,0) \} \varsubsetneqq k^*$. This implies $(1)$. 
\newline 
$(2)$ is clear. We omit the proof. 
\newline 
$(3)$ If $(0,\beta) \in (k^* \times k_{\mathbb{C}}^*)_{reg}$, 
then span$\{ \iota^* u_i ~|~ z_i w_i \ne 0 \} = k^*$ for any $(z,w) \in \mu_{\mathbb{C}}^{-1}(\beta)$. 
By Lemma \ref{criterion}, we have $\mu_{\mathbb{C}}^{-1}(\beta)^{\alpha -ss}=\mu_{\mathbb{C}}^{-1}(\beta)$. 
So we have proved $(3)$.
\newline 
$(4)$ Suppose the contrary. Namely, there exists $(z,w) \in \mu_{\mathbb{C}}^{-1}(\beta)^{\alpha -ss}$ such that 
$(z,w)K_{\mathbb{C}}$ is not closed in $\mu_{\mathbb{C}}^{-1}(\beta)^{\alpha -ss}$. Then there exists 
$(z^{\prime},w^{\prime}) \in \overline{(z,w)K_{\mathbb{C}}} \cap \mu_{\mathbb{C}}^{-1}(\beta)^{\alpha -ss}$ 
such that $(z^{\prime},w^{\prime})K_{\mathbb{C}}$ is closed in $\mu_{\mathbb{C}}^{-1}(\beta)^{\alpha -ss}$. 
By Lemma \ref{moment} we may assume $(z^{\prime},w^{\prime}) \in \mu^{-1}(\alpha,\beta)$.
Since $\dim K_{(z^{\prime},w^{\prime})} \ge 1$, we can easily see that $(z^{\prime},w^{\prime})$ is a critical point of $\mu$. 
Therefore $(\alpha, \beta)$ is a critical value of $\mu$. This is a contradiction. 
\newline
$(5)$ We denote $X(\alpha, \beta)$ by $X(\alpha_1, \alpha_2, \alpha_3)$, where $(\alpha, \beta)=
(\alpha_1, \alpha_2 + \sqrt{-1}\alpha_3)$. Take $(\alpha_1, \alpha_2, \alpha_3), 
(\alpha_1^{\prime}, \alpha_2^{\prime}, \alpha_3^{\prime}) \in (k^* \times k_{\mathbb{C}}^*)_{reg}$. 
We have to show that $X(\alpha_1, \alpha_2, \alpha_3)$ is diffeomorphic to 
$X(\alpha_1^{\prime}, \alpha_2^{\prime}, \alpha_3^{\prime})$. By Lemma \ref{criterion}, we may assume 
that $(0,0, \alpha_3)$ and $(0,\alpha_2^{\prime},0)$ are in $(k^* \times k_{\mathbb{C}}^*)_{reg}$.  
Then $X(\alpha_1, \alpha_2, \alpha_3)$ is diffeomorphic to $X(\alpha_1^{\prime}, \alpha_2, \alpha_3)$, 
because $(0, \alpha_2 + \sqrt{-1}\alpha_3) \in (k^* \times k_{\mathbb{C}}^*)_{reg}$. By rotating 
$I_1$, $I_2$ and $I_3$ cyclically, we see that $X(\alpha_1^{\prime}, \alpha_2, \alpha_3)$ is diffeomorphic to 
$X(\alpha_1^{\prime}, \alpha_2^{\prime}, \alpha_3)$ and also that 
$X(\alpha_1^{\prime}, \alpha_2^{\prime}, \alpha_3)$ is diffeomorphic to 
$X(\alpha_1^{\prime}, \alpha_2^{\prime}, \alpha_3^{\prime})$. Thus we finish the proof. $\hfill \Box$
\vspace{3mm} \newline
It is easy to see when a toric hyperk\"ahler orbifold $X(\alpha,\beta)$ is a manifold.
\begin{proposition}\label{mfd}
Fix $(\alpha,\beta) \in (k^* \times k^*_{\mathbb{C}})_{reg}$. 
Then the subtorus $K$ acts on $\mu^{-1}(\alpha,\beta)$ freely if and only if 
\begin{equation}\label{cmfd}
t^N_{\mathbb{Z}} = k_{\mathbb{Z}} + \sum_{i \not\in J} \mathbb{Z}X_j ~~~ \text{as a $\mathbb{Z}$-module}
\end{equation}
holds for any $J \subset \{ 1, \dots, N \}$ such that $\{ \iota^* u_i ~\vert~ i \in J \}$ forms a basis of $k^*$.
\end{proposition}
See \cite{Ko1} for the proof. Note that (\ref{cmfd}) is the condition for the subtorus $K$. 
It is independent of $(\alpha,\beta)$.

\section{Topology}

In this section we discuss the topology of toric hyperk\"ahler orbifolds, based on the works \cite{BD,Ko1,HS}.
It describes similarities and differences between the topology of ordinary toric orbifolds and 
toric hyperk\"ahler orbifolds very well. 

\subsection{The Core}\label{the core}

Let $\Theta$ be the set of all maps from $\{1, \dots, N\}$ to $\{1,-1 \}$.
For each $\epsilon \in \Theta$, we define 
\begin{equation*}
V_{\epsilon}= \{(z,w) \in \mathbb{H}^N ~|~  
\text{$w_i=0$ if $\epsilon(i)=1$, and $z_i=0$ if $\epsilon(i)=-1$ for $i=1, \dots, N$} \}. 
\end{equation*}
Note that $V_{\epsilon} \cong \mathbb{C}^N$ and $V_{\epsilon} \subset \mu_{\mathbb{C}}^{-1}(0)$. 
If we fix $\alpha \in k^*$, then 
\begin{equation*}
M_{\epsilon}(\alpha)= \{ V_{\epsilon} \cap \mu_{1}^{-1}(\alpha) \}/K \subset X(\alpha,0)
\end{equation*}
is an ordinary toric variety. Moreover, we set 
\begin{equation*}
\Theta_{cpt}(\alpha)= \{ \epsilon \in \Theta ~|~ \text{$M_{\epsilon}(\alpha)$ is compact}  \}.
\end{equation*}
Then we have the following.
\begin{proposition}\label{core}
Fix $\alpha \in k^*$. 
Let $\pi_1 \colon X(\alpha,0) \to X(0, 0)$ be the natural morphism. 
Then the following holds.
\newline
$(1)$ $\pi_1^{-1}([0,0])$ is a $T^n$-equivariant deformation retract of $X(\alpha,0)$. 
\newline
$(2)$ In addition, suppose that $\iota^*u_i \ne 0$ for any $i=1, \dots, N$. 
Then we have 
$$\pi_1^{-1}([0,0]) = \bigcup_{\epsilon \in \Theta_{cpt}(\alpha)}M_{\epsilon}(\alpha).$$ 
\end{proposition}
\textit{Proof.} $(1)$ Define a $\mathbb{C}^{\times}$-action on $X(\alpha,0)$ by $[z,w]s=[sz,sw]$  
for $s \in \mathbb{C}^{\times}$ and $(z,w) \in \mu_{\mathbb{C}}^{-1}(0)^{\alpha -ss}$. 
We also define a $\mathbb{C}^{\times}$-action on $X(0,0)$ in a similar way. 
Then the natural morphism $\pi_1 \colon (X(\alpha,0),I_1) \to (X(0, 0),I_1)$ is 
$(\mathbb{C}^{\times} \times T^n)$-equivariant. So, if we define a 1-parameter 
transformation group $\{ \widetilde{\phi}_t \}$ of $X(\alpha,0)$ by $\widetilde{\phi}_t([z,w])= [e^{-t}z, e^{-t}w]$, then the map $\widetilde{\phi}_t \colon X(\alpha,0) \to X(\alpha,0)$ is $T^n$-equivariant. 
Similarly, we can also define a 1-parameter transformation group $\{ \phi_t \}$ of $X(0,0)$, which is $T^n$-equivariant. 
Since $[0,0]$ is the unique fixed point in $X(0,0)$ for the $\mathbb{C}^{\times}$-action, we have 
$$
\lim_{t \to \infty}\pi_1(\widetilde{\phi}_t([z,w]))=\lim_{t \to \infty}\phi_t(\pi_1([z,w]))=[0,0]
$$
for any $[z,w] \in X(\alpha,0)$.  
Thus we see that $\pi_1^{-1}([0,0])$ is a $T^n$-equivariant deformation retract of $X(\alpha,0)$. 
\newline
$(2)$ Suppose that $\pi_1([z,w])=[0,0]$, that is, 
$$
(0,0) \in \overline{(z,w)K_{\mathbb{C}}}=\overline{\{(z \zeta,w \zeta^{-1}) ~|~ \zeta \in K_{\mathbb{C}}  \}}.
$$
Therefore, we have $z_iw_i=0$ for $i=1, \dots, N$. 
We also have $X=\sum_{i=1}^N a_i X_i \in k$ such that $a_i>0$ if $z_i \ne 0$, and $a_i<0$ if $w_i \ne 0$. 
Since $\iota^*u_j \ne 0$ for any $j=1, \dots, N$, we may assume that $a_j \ne 0$ for any $j=1, \dots, N$. 
If we define $\epsilon \in \Theta$ by $\epsilon (i)=1$ if $a_i>0$, and $\epsilon (i)=-1$ if $a_i<0$, 
then we have $(z,w) \in V_{\epsilon}$ and $\pi_1(M_{\epsilon})=[0,0]$. By Lemma \ref{proper}, $M_{\epsilon}(\alpha)$ is compact. 
So we have $\pi_1^{-1}([0,0]) \subset \bigcup_{\epsilon \in \Theta_{cpt}(\alpha)}
M_{\epsilon}(\alpha)$. By a similar argument 
we have $\pi_1^{-1}([0,0]) \supset \bigcup_{\epsilon \in \Theta_{cpt}(\alpha)}M_{\epsilon}(\alpha)$. $\hfill \Box$
\vspace{3mm} \newline 
\begin{remark}
\rm{If $\iota^* u_N=0$, then we have $K \subset T^{N-1}$ and $X(\alpha,0)$ is decomposed into 
the product of a lower dimensional toric hyperk\"ahler variety and $\mathbb{H}$. Therefore the 
assumption of Proposition \ref{core} can be dropped easily if we slightly change the statement. 
See \cite{BD, HS} and the next subsection for detail.}
\end{remark}
By Proposition \ref{core} we call the set $\pi_1^{-1}([0,0])$ the {\it core} of the toric hyperk\"ahler orbifold $X(\alpha,0)$.
The above proof of $(1)$ is based on the argument in \cite {Kr,N1}, where the similar result is proved in the case of quiver varieties. 
Bielawski and Dancer proved $(1)$ and $(2)$ in a slightly different formulation \cite{BD}.
\vspace{3mm}
\newline
\textit{Proof of Proposition \ref{proper}}. 
Let us prove the surjectivity of the natural morphism $\pi_1 \colon X(\alpha,\beta) \to X(0, \beta)$. 
Take an arbitrary point $(a,b) \in \mu^{-1}(0,\beta)$. If we set 
$W=\sum_{i=1}^N \mathbb{R}|a_i|^2 \iota^* u_i + \sum_{i=1}^N \mathbb{R}|b_i|^2 (-\iota^* u_i) \subset k^*$, then we have 
\begin{equation}\label{W}
W=\{ v \in k^* ~|~ \langle v, \underline{k}_{(a,b)} \rangle =0 \}, 
\end{equation}
where $\underline{k}_{(a,b)}=Lie K_{(a,b)} \subset k$ is the Lie algebra of the isotropy subgroup $K_{(a,b)} \subset K$ at $(a,b)$. 
Set $J= \{ i ~|~\iota^* u_i \not\in W \}$. Since $a_i=b_i=0$ for any $i \in J$, by Lemma \ref{criterion}, 
we have $0 \in \sum_{i \not\in J} \mathbb{R}_{>0}|a_i|^2 \iota^* u_i + \sum_{i \not\in J} \mathbb{R}_{>0}|b_i|^2 (-\iota^* u_i)$.
Therefore we have 
\begin{equation}\label{W2}
W= \sum_{i \not\in J} \mathbb{R}_{>0}|a_i|^2 \iota^* u_i + \sum_{i \not\in J} \mathbb{R}_{>0}|b_i|^2 (-\iota^* u_i).
\end{equation}

On the other hand, $\alpha \in k^*$ induces $\underline{\alpha} \in (\underline{k}_{(a,b)})^*$. 
Let $\underline{X}(\underline{\alpha},\underline{0})$ and $\underline{X}(\underline{0},\underline{0})$ be hyperk\"ahler quotients of $\mathbb{H}^J=\{ (z,w) \in \mathbb{H}^N~|~ \text{$z_i=w_i=0$ for any $i \not\in J$} \}$ by $K_{(a,b)}$ at $(\underline{\alpha},\underline{0}),(\underline{0},\underline{0}) \in (\underline{k}_{(a,b)})^* \times (\underline{k}_{(a,b)\mathbb{C}})^*$ respectively. 
Let $\underline{\pi}_1 \colon \underline{X}(\underline{\alpha},\underline{0}) \to \underline{X}(\underline{0},\underline{0})$ be the natural morphism. By the proof of Proposition \ref{core} (2)
we can take a point $[c,d] \in (\underline{\pi}_1)^{-1}([0,0])$, where $(c,d) \in \mathbb{H}^J$. 
Namely, we have $(0,0) \in \overline{(c,d)K_{(a,b)}} $. If we set $(z,w)=(a,b)+(c,d) \in \mathbb{H}^N$, 
then we have $(a,b) \in \overline{(z,w)K_{(a,b)}}$. 
It is obvious that $(z,w) \in \mu_{\mathbb{C}}^{-1}(\beta)$. 
Moreover, we have 
\begin{equation}\label{W3}
\underline{\alpha} \in \sum_{i \in J} \mathbb{R}_{\ge 0}|c_i|^2 \underline{\iota}^* u_i + \sum_{i \in J} \mathbb{R}_{\ge 0}|d_i|^2 (-\underline{\iota}^* u_i),
\end{equation}
where $\underline{\iota} \colon \underline{k}_{(a,b)} \to t^{J}$ is the natural inclusion and $t^J=\mathrm{span}\{ X_i ~|~ i \in J \} \subset t^N$.  
Therefore, by (\ref{W}), (\ref{W2}) and (\ref{W3}), we have $(z,w) \in \mu_{\mathbb{C}}^{-1}(\beta)^{\alpha -ss}$.
Thus we proved $\pi_1([z,w])=[a,b]$.  $\hfill \Box$
\vspace{3mm} \newline 

\subsection{The Betti Numbers and the Cohomology Rings}

 By Proposition \ref{core}, to study the topology of a toric hyperk\"ahler orbifold $X(\alpha,0)$, it is 
enough to study the topology of its core.  
To do this, the action of the quotient torus plays an important role. 
Actually Bielawski and Dancer computed the Betti numbers  \cite{BD}, 
and the author determined the cohomology ring \cite{Ko1} by using it. 
In this section we review these results. 

Fix $(\alpha,0) \in (k^* \times k^*_{\mathbb{C}})_{reg}$. 
First we describe the hyperk\"ahler moment map 
$$
\mu_{T^n}=(\mu_{T^n,1},\mu_{T^n \mathbb{C}}) \colon X(\alpha,0) 
\to (t^n)^* \otimes \mathbb{R}^3 \cong (t^n)^* \times (t^n_{\mathbb{C}})^*
$$ 
for the action of the quotient torus $T^n =T^N/K$ on $X(\alpha,0)$. 
Fix $h \in (t^N)^*$ such that $\iota^* h= \alpha$. 
Then the hyperk\"ahler moment map $\mu_{T^n} \colon X(\alpha,0) 
\to (t^n)^* \otimes \mathbb{R}^3$ is given by
\begin{align*}
& \mu_{T^n,1}([z,w]) = \pi\sum_{i=1}^N (|z_i|^2-|w_i|^2)u_i-h \in (t^n)^* \subset (t^N)^*, \\
& \mu_{T^n \mathbb{C}}([z,w]) 
 = -2 \pi \sqrt{-1} \sum_{i=1}^N z_i w_i u_i \in (t^n_{\mathbb{C}})^* \subset (t^N_{\mathbb{C}})^*.
\end{align*}
Note that $\mu_{T^N,1}(z,w)-h \in \mathrm{Ker} \{ \iota^* \colon (t^N)^* \to k^* \} = \mathrm{Im} 
\{ \pi^* \colon (t^n)^* \to (t^N)^* \}$. Similarly, we have $-2 \pi \sqrt{-1} \sum_{i=1}^N z_i w_i u_i \in 
\mathrm{Ker}~ \iota^*  \otimes \mathbb{C}=\mathrm{Im}~\pi^* \otimes \mathbb{C}$.

We define a hyperplane $F_i \subset (t^n)^*$ by 
\[
F_i= \{ p \in (t^n)^* \vert \langle  \pi^* p +h, X_i \rangle =0 \}
\quad \text{for $i=1, \dots, N$.}
\]
We note that $F_i = \emptyset$ if $\pi(X_i)=0$, because we assume $(\alpha,0) \in (k^* \times k^*_{\mathbb{C}})_{reg}$. 
Let $\Theta$ be the set of all maps from $\{1, \dots, N\}$ to $\{1,-1 \}$ as in the last section. 
Then these hyperplanes divide $(t^n)^*$ into a finite number of polyhedra 
$\{ \Delta_{\epsilon} \vert \epsilon \in \Theta \}$, where  
$\Delta_{\epsilon} \in (t^n)^*$ is defined by
\[
\Delta_{\epsilon}=\{ p \in (t^n)^* \vert 
\epsilon(i) \langle \pi^* p + h, X_i \rangle   \ge 0  \hspace{3mm} \text{for any $i =1, \dots, N$} \}.
\]
We note that some $\Delta_{\epsilon}$ may be empty. 
Then we have the following lemma.
\begin{lemma}
$[z,w] \in \mu_{T^n}^{-1}(\Delta_{\epsilon},0)$  if and only if $(z,w) \in V_{\epsilon} 
\cap \mu_{1}^{-1}(\alpha)$. 
Namely, $M_{\epsilon}(\alpha)=\mu_{T^n}^{-1}(\Delta_{\epsilon},0)$ holds. 
\end{lemma}
\textit{Proof.}  First we note that $[z,w] \in \mu_{T^n}^{-1}((t^n)^*,0)$ if and only if 
$z_i w_i=0$ for all $i=1, \dots, N$. If we set $p=\mu_{T^n,1}([z,w])$, we have 
\[
\langle \pi^* p +h,X_i \rangle = \pi (\vert z_i \vert^2 - \vert w_i \vert^2 ) \quad \text{for $i=1, \dots, N$.}
\]
Then the lemma follows. $\hfill \Box$
\vspace{3mm} \newline 

We want to study the topology of the core. 
Note that each component of the core is an ordinary toric orbifold, even if we do not assume that $\iota^* u_i \ne 0$ for any $i=1, \dots, N$. 
So we have to investigate how these components intersect to each other. 
This is described by the polyhedral complex, which is define by 
$$
\mathcal{C}(X(\alpha,0))= \{ \sigma ~|~  \text{$\sigma $ is a {\it compact} face of a polyhedron $\Delta_{\epsilon}$ for some $\epsilon \in \Theta $} \}.
$$
It should be remarked that, to define $\mathcal{C}(X(\alpha,0))$, we need $h \in (t^N)^*$ 
such that $\iota^* h=\alpha$. However, $\mathcal{C}(X(\alpha,0))$ is determined by $\alpha$ up to 
parallel translation. So we use this notation. 

Note that the polyhedral complex $\mathcal{C}(X(\alpha,0))$ is determined by the hyperplanes $F_1, \dots, F_N$.  
Therefore, the cohomology of the toric hyperk\"ahler orbifold $X(\alpha,0)$ is described in terms of 
the hyperplanes $F_1, \dots, F_N$ or the polyhedral complex $\mathcal{C}(X(\alpha,0))$. 
In fact, we have the following. 

\begin{theorem}\label{betti}{\rm \cite{BD}}
Fix $(\alpha,0) \in (k^* \times k_{\mathbb{C}}^*)_{reg}$. 
Let $d_k$ be the number of $k$-simplexes of the polyhedral complex 
$\mathcal{C}(X(\alpha,0))$. Then the Poincar\'{e} polynomial $P_t(X(\alpha,0))$ is given by 
$$
P_t(X(\alpha,0))= \sum_{k=0}^n d_k (t^2 -1)^k.
$$
\end{theorem}
Note that the Poincar\'{e} polynomial of an ordinary toric orbifold is the same form if we replace 
$\mathcal{C}(X(\alpha,0))$ by the convex polytope associated to the toric orbifold. 
See Section 4.5 in \cite{Ful}.

In Proposition \ref{mfd} we stated the condition so that $X(\alpha, \beta)$ is a smooth manifold 
for $(\alpha,\beta) \in (k^* \times k_{\mathbb{C}}^*)_{reg}$. 
Recall the condition is the one for the subtorus $K$, and independent of  
$(\alpha,\beta)$. Under the condition $\mu^{-1}(\alpha, \beta)$ can be viewed as a 
principal $K$-bundle on $X(\alpha, \beta)$. Let $L_i$ be a line bundle on 
$X(\alpha, \beta)$ associated the character $\iota^* u_i$ of $K$ for $i=1, \dots, N$. 

\begin{lemma}\label{str}
Suppose that $\Delta_{\epsilon} \cap F_i$ is a face of $\Delta_{\epsilon}$ of codimension one. 
Then the homology class represented by $\mu_{T^n}^{-1}(\Delta_{\epsilon} \cap F_i,0)$ is 
the Poincar\'{e} dual of $\epsilon (i)c_1 (L_{i})$ in $M_{\epsilon}(\alpha)$. 
\end{lemma}
\textit{Proof.} Suppose $[z,w] \in \mu_{T^n}^{-1}((t^n)^*,0)$. Then $[z,w] \in F_i$ 
if and only if $z_i = w_i =0$. Then we have
\[
\mu_{T^n}^{-1}(\Delta_{\epsilon} \cap F_i,0)= 
\begin{cases}
\{ [z,w] \in M_{\epsilon}(\alpha) ~|~ z_i=0 \} \quad 
\text{if $\epsilon (i)=1$,} \\
\{ [z,w] \in M_{\epsilon}(\alpha) ~|~ w_i=0 \} \quad 
\text{if $\epsilon (i)=-1$.}
\end{cases}
\]
Note that the divisor defined by the equation $z_i=0$ in $(X(\alpha,\beta),I_1)$ corresponds to $L_{i}$ and that the divisor defined by the equation $w_i=0$ in $(X(\alpha,\beta),I_1)$ corresponds to $L_{i}^*$. 
So we finish the proof.  $\hfill \Box$
\vspace{3mm} \newline 

The cohomology ring of a toric hyperk\"ahler manifold $X(\alpha,\beta)$ is described as follows. 

\begin{theorem}\label{ring}{\rm \cite{Ko1}}
Suppose that $(\alpha,\beta) \in (k^* \times k_{\mathbb{C}}^*)_{reg}$ 
and $X(\alpha, \beta)$ is a smooth manifold.
Define a ring homomorphism $\Phi \colon \mathbb{Z}[u_1, \dots, u_N] \to 
H^*(X(\alpha,\beta); \mathbb{Z})$ by $\Phi(u_i)=c_1(L_i)$. 
Then the following holds:
\newline
$(1)$ The map $\Phi$ is surjective. So we have a ring isomorphism
\[
H^*(X(\alpha,\beta);\mathbb{Z}) \cong 
\mathbb{Z}[u_1, \dots, u_N]/ \ker \Phi.
\]
$(2)$ The ideal $\ker \Phi$ is generated by the following two types of elements:
\newline
\hspace*{5mm} \textnormal{$(i)$~(linear relations)}~~ $\displaystyle \sum_{i=1}^N a_i u_i \in \ker \iota^* \cap 
(t^N_{\mathbb{Z}})^*$, 
\newline
\hspace*{5mm} \textnormal{$(ii)$~(nonlinear relations)}~~ $\displaystyle \prod_{i \in J}u_i$~~ for 
all $J \subset \{ 1, \dots, N \}$ such that  $\bigcap_{i \in J}F_i = \emptyset$.
\end{theorem}
\begin{remark}
\rm{For $\epsilon \in \Theta_{cpt}(\alpha)$ we have $H(M_{\epsilon}(\alpha),\mathbb{Z}) \cong 
\mathbb{Z}[u_1, \dots, u_N]/I_{\epsilon}$, where the ideal $I_{\epsilon}$ are generated by the following 
two types of elements;
\newline
\hspace*{5mm} $(a)$ the same as $(i)$ in Theorem \ref{ring} $(2)$, 
\newline
\hspace*{5mm} $(b)$ $\displaystyle \prod_{i \in J}u_i$~~ for 
all $J \subset \{ 1, \dots, N \}$ such that 
$\displaystyle (\bigcap_{i \in J}F_i) \cap \Delta_{\epsilon}= \emptyset$. 
\newline
To compute the cohomology of the core, we have to take all $M_{\epsilon}(\alpha)$ for 
$\epsilon \in \Theta_{cpt}(\alpha)$ into account at same time. This is the geometric meaning of 
the nonlinear relations in Theorem \ref{ring}. }
\end{remark}
\begin{remark} \rm{In \cite{HS} Hausel and Sturmfels generalized Theorem \ref{ring} to orbifold cases, where the coefficient $\mathbb{Z}$ is replaced by $\mathbb{R}$. In Subsection \ref{other} we discuss their proof. 
Our proof, which will be described in the next subsection, also works for orbifold cases 
if $L_i$ is considered as an orbi-line bundle on $X(\alpha, \beta)$. }
\end{remark}

\subsection{Outline of the proofs of Theorems \ref{betti} and \ref{ring}}
Let us sketch the idea of the proofs of Theorems \ref{betti} and \ref{ring}. 
By Proposition \ref{indep} it is enough to compute the cohomology ring of $X(\alpha,0)$ for 
$(\alpha,0) \in (k^* \times k^*_{\mathbb{C}})_{reg}$.
Let $F_1, \dots, F_N$ be the hyperplanes associated to $X(\alpha,0)$. 
Let $Y_1$ and $\mathcal{C}(Y_1)$ be another toric hyperk\"ahler variety and the polyhedral complex 
associated to the hyperplanes $F_1, \dots, F_{N-1}$ respectively.  
By Proposition \ref{indep}, the topology of $X(\alpha,0)$ does not change as far as $(\alpha,0) \in 
(k^* \times k_{\mathbb{C}}^*)_{reg}$. Since the variation of $\alpha$ corresponds to parallel 
translations of the hyperplanes, we can choose $\alpha$ so that $\mathcal{C}(Y_1)$ is contained 
in one of the components of $k\* \setminus F_N$. Then we consider another toric hyperk\"ahler variety 
$Y_2$ associated to the hyperplanes $F_1 \cap F_N, \dots, F_{N-1}\cap F_N$. 
Bielawski and Dancer expressed the Betti numbers of $X(\alpha,0)$ in terms of those of $Y_1$ and $Y_2$ by the Meyer-Vietoris argument. 
Then by an inductive argument they computed the Betti numbers of $X(\alpha,0)$ \cite{BD}. 
In \cite{Ko1} the author determined the cohomology ring of $X(\alpha,0)$ 
from the one of $Y_1$ and $Y_2$ by a refinement of the above argument. 

\subsection{Another proof}\label{other}
Another proof was given by Hausel and Sturmfels. In this subsection we explain their proof. 
Let us consider a symplectic quotient $L(\alpha)=\mu_1^{-1}(\alpha)/K$ of $\mathbb{H}^N$ by $K$, 
which is a non-compact ordinary toric variety. Note that a toric hyperk\"ahler variety 
$(X(\alpha,\beta),I_1)=(\mu_1^{-1}(\alpha) \cap \mu_{\mathbb{C}}^{-1}(\beta))/K$ is a subvariety of $L(\alpha)$. 
Hausel and Sturmfels defined the core of $L(\alpha)$, which is a deformation retract 
of $L(\alpha)$ for $(\alpha,0) \in (k^* \times k^*_{\mathbb{C}})_{reg}$. Moreover, they observed 
that the core of $X(\alpha,0)$ is the same as the one of $L(\alpha)$. 
So they proved that the cohomology ring of $X(\alpha,0)$ is isomorphic to the one of $L(\alpha)$, which 
can be computed by the standard argument for ordinary toric varieties \cite{Ful}. 
The argument can be considered as a refinement of the Morse theory. So the properness of a Morse function is 
important. Since $L(\alpha)$ is non-compact, Hausel and Sturmfels carefully choose an appropriate 
Morse function and computed the cohomology ring of $X(\alpha,0)$ \cite{HS}. 

Here we prove that $X(\alpha,0)$ is homotopy equivalent to $L(\alpha)$ in a different way from them. 

\begin{proposition}
Suppose $(\alpha,0) \in (k^* \times k^*_{\mathbb{C}})_{reg}$. Then $L(\alpha)$ is diffeomorphic to 
the product $X(\alpha,0) \times k^*_{\mathbb{C}}$. In particular, $X(\alpha,0)$ is homotopy equivalent to 
$L(\alpha)$.
\end{proposition}
\textit{Proof.} Since $(\alpha,0) \in (k^* \times k^*_{\mathbb{C}})_{reg}$, due to Lemma \ref{indep}, the natural morphism 
$\pi_2 \colon (X(\alpha, \alpha_2,\alpha_3),I_2) \to (X(\alpha, 0,\alpha_3),I_2)$ is a biholomorphic map for any $\alpha_2,\alpha_3 \in k^*$. 
Similarly, $\pi_3 \colon (X(\alpha, 0,\alpha_3),I_3) \to (X(\alpha, 0, 0),I_3)$ is also a biholomorphic map. 
So we have a diffeomorphism $\pi_3 \circ \pi_2 \colon X(\alpha,\beta) \to X(\alpha,0)$ for any 
$\beta \in k^*_{\mathbb{C}}$. Since $L(\alpha)$ is equal to $\bigcup_{\beta \in k^*_{\mathbb{C}}}X(\alpha,\beta)$, 
this map induces a diffeomorphism from $L(\alpha)$ to $X(\alpha,0) \times k^*_{\mathbb{C}}$. $\hfill \Box$
\vspace{3mm} \newline 

We refer the reader to \cite{P} for other methods to compute the cohomology ring and also for the computation of the equivariant cohomology of $X(\alpha,\beta)$ for the torus action. 

\subsection{Another presentation of the cohomology ring}
The topology of a toric hyperk\"ahler variety $X(\alpha,\beta)$ for $(\alpha,\beta) \in (k^* \times k^*_{\mathbb{C}})_{reg}$ is determined only by the subtorus $K$. 
So we have a presentation of its cohomology ring in terms of $K$ as follows. 

\begin{theorem}\label{ring2}
Suppose the same assumption as in Theorem \ref{ring}. Let $\Lambda^{(1)}= \{ H^{(1)}_1, \dots, H^{(1)}_l \}$ 
be the set of 1-dimensional isotropy subgroups in $K$. Fix a non-zero element $Y_s \in Lie H_s^{(1)}$ and 
set $J_s = \{ i \in \{ 1, \dots, N \} ~|~ \langle u_i, Y_s \rangle \ne 0 \}$ for $s=1, \dots, l$.  
Then \textnormal{$(ii)$~(nonlinear relations)} in Theorem \ref{ring} $(2)$ can be replaced by the following;
\newline
\hspace*{5mm} \textnormal{$(ii)^{\prime}$~(nonlinear relations)}~~ $\displaystyle \prod_{i \in J_s}u_i$ ~~for 
$s=1, \dots, l$. 
\end{theorem}
Since the above presentation is slightly different from the one in \cite{Ko1}, we give a proof here, assuming 
Theorem \ref{ring}. 
\vspace{3mm} \newline \textit{Proof.} Let $I$ and $I^{\prime}$ be the ideals in $\mathbb{Z}[u_1, \dots, u_N]$ in Theorem \ref{ring} 
and Theorem \ref{ring2} respectively. First we show $I^{\prime} \subset I$. We have to show 
$\bigcap_{i \in J_s} F_i = \emptyset$ for any $s=1, \dots, l$. 
Suppose that there exists $p \in \bigcap_{i \in J_s} F_i$, that is, 
$\langle \pi^* p+h,X_i \rangle =0$ for all $i \in J_s$. Therefore we have
$\pi^* p+h \in \mathrm{span} \{u_i ~|~ i \not\in J_s \}$. Thus we have 
$\alpha \in \mathrm{span} \{\iota^* u_i ~|~ i \not\in J_s \}$. 
Since $(\alpha,0) \in (k^* \times k^*_{\mathbb{C}})_{reg}$, 
according to Lemma \ref{indep}, we have $\mathrm{span} \{\iota^* u_i ~|~ i \not\in J_s \}=k^*$. 
Since $\langle Y_s,~ \iota^* u_i \rangle =0$ for all $i \not\in J_s$, we have $Y_s=0$. 
This is a contradiction. Thus we see that $\bigcap_{i \in J_s} F_i = \emptyset$. 

Next we show $I \subset I^{\prime}$. Suppose that $\bigcap_{i \in J}F_i= \emptyset$. 
We have to prove that there exists $s \in \{1, \dots, l \}$ such that $J_s \subset J$. 
First we prove that $\mathrm{span} \{ \iota^* u_i ~|~ i \not\in J \} \subsetneqq k^*$. 
Suppose that 
$\mathrm{span} \{ \iota^* u_i ~|~ i \not\in J \} = k^*$. 
Then we can choose $A \subset J^c$, where $J^c = \{ 1, \dots, N \} \setminus J$, such that $\{ \iota^* u_i ~|~ i \in A \}$ forms a basis of $k^*$. 
So we can write $\alpha = \sum_{i \in A}c_i \iota^*u_i$. 
Therefore there exists $(z,w) \in \mu^{-1}(\alpha,0)$ such that $z_i =w_i =0$ if $i \not\in A$, 
$w_i =0$ if $i \in A$ and $c_i \ge 0$, $z_i =0$ if $i \in A$ and $c_i \le 0$. 
If we set $p=\mu_{T^n,1}([z,w])$, according to Lemma \ref{str}, we have 
$p \in \bigcap_{i \in A^c} F_i \subset \bigcap_{i \in J} F_i$. This is a contradiction. 
Thus we proved $\mathrm{span} \{ \iota^* u_i ~|~ i \not\in J \} \subsetneqq k^*$. 
So there exists a wall $W_s$ such that $\mathrm{span} \{ \iota^* u_i ~|~ i \not\in J \} \subset W_s$. 
This implies $J_s \subset J$. $\hfill \Box$
\vspace{3mm} \newline 
\subsection{The hyperk\"ahler Kirwan map}
Assume that a toric hyperk\"ahler variety $X(\alpha,\beta)$ is smooth. 
Let $H_K^*(\mathbb{H}^N; \mathbb{Z})$ be the $K$-equivariant cohomology of $\mathbb{H}^N$. 
Define the map 
$$
\kappa_{\mathbb{Z}} \colon H_K^*(\mathbb{H}^N; \mathbb{Z}) \to H^*(X(\alpha,\beta); \mathbb{Z})
$$ 
to be the composition of the restriction map 
$r \colon H_K^*(\mathbb{H}^N; \mathbb{Z}) \to H_K^*(\mu^{-1}(\alpha,\beta); \mathbb{Z})$ 
and the natural isomorphism $i \colon H_K^*(\mu^{-1}(\alpha,\beta); \mathbb{Z}) \to 
H^*(X(\alpha,\beta); \mathbb{Z})$. 
The map $\kappa_{\mathbb{Z}}$ is an hyperk\"ahler analogue of the Kirwan map for a 
symplectic quotient. Kirwan proved that the Kirwan map is surjective for symplectic 
quotients in $\mathbb{R}$-coefficients (even if the quotient is smooth) under very weak assumptions. 
Theorem \ref{ring} implies that the hyperk\"ahler Kirwan map is surjective for toric hyperk\"ahler manifolds 
even in $\mathbb{Z}$-coefficients, because $H_K^*(\mathbb{H}^N; \mathbb{Z}) \cong S^* (k_{\mathbb{Z}}^*) \cong
\mathbb{Z}[u_1, \dots, u_N]/ (\ker \iota^* \cap (t^N_{\mathbb{Z}})^*)$, 
where $S^* (k_{\mathbb{Z}}^*)$ is the symmetric power of $k_{\mathbb{Z}}^*$. 
In the orbifold case the hyperk\"ahler Kirwan map 
$$
\kappa \colon H_K^*(\mathbb{H}^N; \mathbb{R}) \to H^*(X(\alpha,\beta); \mathbb{R})
$$
is also surjective. 
In \cite{Ko2} the author computed the cohomology ring and proved the surjectivity 
of the hyperk\"ahler Kirwan map for the hyperk\"ahler polygon spaces, which are hyperk\"ahler quotients 
of $\mathbb{H}^N$ by certain non-abelian Lie groups. 
It is not known whether the hyperk\"ahler Kirwan map is surjective or not in general. 
Surjectivity of the hyperk\"ahler Kirwan map is one of the central topics for study of the topology of hyperk\"ahler quotients.

\section{Examples}
In this section we give two examples, which illustrate the results in other sections. The first one 
is the most fundamental. 
\begin{example}\label{ex1} 
\rm{Consider the subtorus $K$ of $T^{n+1}$ whose Lie algebra $k$ is spanned by $X_1 + \dots +X_{n+1}$. 
If we denote the dual basis by $v \in k^*$, then we have 
$\iota^* u_1 = \dots = \iota^* u_{n+1} =v$.  

The hyperk\"ahler moment map for the action of $K$ on $\mathbb{H}^{n+1}$
$$
\mu=(\mu_1, \mu_{\mathbb{C}}) \colon \mathbb{H}^{n+1} \to k^* \otimes 
\mathbb{R}^3 \cong k^* \times k^*_{\mathbb{C}}
$$
is given by 
\begin{align*}
\mu_1(z,w)= & \pi \sum_{i=1}^{n+1} (\vert z_i \vert^2 - \vert w_i \vert^2)v \in k^*, \\
\mu_{\mathbb{C}}(z,w)= & -2\pi \sqrt{-1} (\sum_{i=1}^{n+1} z_i w_i) v \in k^*_{\mathbb{C}}. 
\end{align*}
Then we have $(k^* \times k^*_{\mathbb{C}})_{reg}=k^* \times k^*_{\mathbb{C}} \setminus \{(0,0) \}$.

First we consider the case $\beta=0$. Fix $\alpha_+ =a_+ v$ and $\alpha_- =a_- v$, where $a_+ >0, ~a_- <0$. 
By Lemma \ref{criterion} we have \begin{equation}\label{plus}
\mu_{\mathbb{C}}^{-1}(0)^{\alpha_+ -ss}= \{ (z,w) \in \mu_{\mathbb{C}}^{-1}(0) ~|~ z \ne 0 \}. 
\end{equation}
By Proposition \ref{indep} we have $(X(\alpha_+, 0),I_1)
=\mu_{\mathbb{C}}^{-1}(0)^{\alpha_+ -ss} / K_{\mathbb{C}}$, which is biholomorphic to the total space of the 
cotangent bundle of $\mathbb{C}P^n$.
Similarly, we have 
\begin{equation}\label{minus}
\mu_{\mathbb{C}}^{-1}(0)^{\alpha_- -ss}=\{ (z,w) \in \mu_{\mathbb{C}}^{-1}(0) ~|~ w \ne 0 \}. 
\end{equation}
Therefore we have $(X(\alpha_-, 0), I_1)=\mu_{\mathbb{C}}^{-1}(0)^{\alpha_- -ss} / K_{\mathbb{C}}$, 
which is also biholomorphic to $T^*\mathbb{C}P^n$.

Next we consider the case $\beta \ne 0$. Fix an arbitrary $\alpha$. Then, by Proposition \ref{indep}, 
$(X(\alpha,\beta),I_1) \cong (X(0,\beta),I_1)$ is an affine variety. 
Therefore $(X(\alpha,\beta),I_1)$ is diffeomorphic to $T^* \mathbb{C}P^n$, but not biholomorphic as a 
complex manifold. 

By Theorems \ref{ring} and \ref{ring2} we have  
$H^*(X(\alpha,\beta),\mathbb{Z}) \cong \mathbb{Z}[u_1, \dots, u_{n+1}]/I$, where the ideal $I$ is generated by 
two types of elements. The linear relations are generated by $u_1-u_2,u_2-u_3, \dots u_n -u_{n+1}$.  
The nonlinear relation is generated by $u_1 u_2 \dots u_{n+1}$. 
Therefore we have 
\begin{align*}
& H^*(X(\alpha,\beta),\mathbb{Z}) \\ 
& \cong \mathbb{Z}[u_1, \dots, {n+1}]/(u_1-u_2,u_2-u_3, \dots u_n -u_{n+1}, ~u_1 u_2 \dots u_{n+1})\\
& \cong \mathbb{Z}[v]/( v^{n+1}).
\end{align*}}
\end{example}
\begin{example}\label{ex2}
\rm{Consider the subtorus $K$ of $T^5$ whose Lie algebra $k$ is spanned by $X_1 +X_2 +X_4$ and $X_1 +X_3 +X_5$. 
Let $\{ v_1,v_2 \}$ be the dual basis. Then we have $\iota^* u_1 =v_1 +v_2$, $\iota^* u_2 =\iota^* u_4 =v_1$ 
and $\iota^* u_3 =\iota^* u_5 =v_2$. 

The hyperk\"ahler moment map for the action of $K$ on $\mathbb{H}^5$
$$
\mu=(\mu_1, \mu_{\mathbb{C}}) \colon \mathbb{H}^5 \to k^* \otimes 
\mathbb{R}^3 \cong k^* \times k^*_{\mathbb{C}}
$$
is given by 
\begin{align*}
\mu_1(z,w)= & \pi \sum_{i=1,2,4} (\vert z_i \vert^2 - \vert w_i \vert^2)v_1 
+ \pi \sum_{i=1,3,5} (\vert z_i \vert^2 - \vert w_i \vert^2)v_2
\in k^*, \\
\mu_{\mathbb{C}}(z,w)= & -2\pi \sqrt{-1} (\sum_{i=1,2,4} z_i w_i) v_1 
-2\pi \sqrt{-1} (\sum_{i=1,3,5} z_i w_i) v_2 \in k^*_{\mathbb{C}}. 
\end{align*}

We have three 1-dimensional isotropy subgroups $H^{(1)}_1, H^{(2)}_1, H^{(3)}_1$, whose Lie algebras are
spanned by $Y_1= X_1 +X_3 +X_5$, $Y_2=X_1 +X_2 +X_4$ and $Y_3= X_2 -X_3 +X_4 -X_5$ respectively.  
So we have three walls $W_1 =\mathbb{R}v_1$, $W_2 =\mathbb{R}v_2$ and $W_3 =\mathbb{R}(v_1 +v_2)$ in $k^*$. 
By Proposition \ref{indep}, we see that $(\alpha,0) \in (k^* \times k^*_{\mathbb{C}})_{reg}$ if and only if 
$\alpha \in k^* \setminus \bigcup_{i=1,2,3}W_i$. Connected components of this set are called chambers. 
So we have six chambers $\mathcal{C}_1, \dots, \mathcal{C}_6$ in $k^*$ as in Figure 1. 
(This is the chamber structure for $\beta=0$. See Section \ref{variation} for the precise definition.)

\setlength{\unitlength}{1mm}
\begin{picture}(100,52)(-50,-27)
%\put(-50,-27){\circle*{1}}
%\put(-50,24){\circle*{1}}
%\put(50,-27){\circle*{1}}
%\put(50,24){\circle*{1}}

\put(-20,0){\line(1,0){40}}
\put(0,-20){\line(0,1){40}}
\put(-20,-20){\line(1,1){40}}
\put(10,0){\circle*{1}}
\put(0,10){\circle*{1}}
\put(10,10){\circle*{1}}
\put(5,-4){{$\iota^* u_2 = \iota^* u_4=v_1$}}
\put(-27,9){{$\iota^* u_3 = \iota^* u_5=v_2$}}
\put(12,9){{$\iota^* u_1=v_1+v_2$}}
\put(13,4){{$\mathcal{C}_1$}}
\put(5,15){{$\mathcal{C}_2$}}
\put(-14,15){{$\mathcal{C}_3$}}
\put(-14,-6){{$\mathcal{C}_4$}}
\put(-7,-14){{$\mathcal{C}_5$}}
\put(13,-14){{$\mathcal{C}_6$}}
\put(-6,-24){{Figure 1}}
\put(21,0){{$W_1$}}
\put(21,19){{$W_3$}}
\put(-5,19){{$W_2$}}
\end{picture}

If we take $\alpha_1=sv_1+tv_2 \in \mathcal{C}_1$, then we have $s > t > 0$. Moreover, we have $\Theta_{cpt}(\alpha_1)= \{\epsilon_1, \epsilon_2 \}$, where 
\begin{equation*}
\begin{split}
\epsilon_1 (i) & =1 \hspace{2mm} \text{for $i=1,2,3,4,5$}
\end{split}, \hspace{5mm}
\begin{split}
\epsilon_2 (i) & =\left\{
\begin{array}{ll}
1 & \hspace{2mm} \text{for $i=1,2,4$} \\
-1 & \hspace{2mm} \text{for $i=3,5$.} \\
\end{array} \right.
\end{split}
\end{equation*}
The associated polyhedral complex $\mathcal{C}(X(\alpha_1,0))$ consists of 
all faces of $\Delta_{\epsilon_1}$ and $\Delta_{\epsilon_2}$ 
as in Figure 2, where we take an appropriate coordinate $(a_1,a_2,a_3)$ in $(t^3)^*$ such that 
$F_i = \{ (a_1,a_2,a_3) ~|~  a_i =0 \} $ for $i=1,2,3$, $F_4 = \{ (a_1,a_2,a_3) ~|~  a_1 + a_2 = s \}$ and 
$F_5 = \{ (a_1,a_2,a_3) ~|~  a_1+a_3 =t \} $.  
Note that $F_1 \cap F_2 \cap F_4= F_1 \cap F_3 \cap F_5= F_2 \cap F_3 \cap F_4 \cap F_5=\emptyset$. 

Let us compute the cohomology ring of $X(\alpha_1,0)$. By Theorems \ref{ring} and \ref{ring2} we have  
$H^*(X(\alpha_1,0),\mathbb{Z}) \cong \mathbb{Z}[u_1, \dots, u_5]/I$, where the ideal $I$ is generated by 
two types of elements. The linear relations are generated by $u_2-u_4$, $u_3-u_5$ and $u_1-u_2-u_3$.  
The nonlinear relations are generated by $u_1 u_2 u_4$, $u_1 u_3 u_5$ and $u_2 u_3 u_4 u_5$. 
Therefore we have 
\begin{align*}
& H^*(X(\alpha_1,0),\mathbb{Z}) \\ 
& \cong \mathbb{Z}[u_1, \dots, u_5]/(u_2-u_4, u_3-u_5, u_1-u_2-u_3, u_1 u_2 u_4, u_1 u_3 u_5, u_2 u_3 u_4 u_5)\\
& \cong \mathbb{Z}[v_1,v_2]/( (v_1+v_2)v_1^2, (v_1+v_2)v_2^2, v_1^2v_2^2).
\end{align*}

\setlength{\unitlength}{1mm}
\begin{picture}(110,50)(0,0)
%\put(0,0){\circle*{1}}
%\put(110,0){\circle*{1}}
%\put(0,50){\circle*{1}}
%\put(110,50){\circle*{1}}

\put(10,10){\line(1,1){30}}
\put(10,20){\line(3,1){30}}
\put(5,15){\line(1,1){15}}
\put(20,30){\line(0,1){15}}
\put(20,40){\line(1,0){20}}
\put(40,30){\line(0,1){10}}
\put(10,10){\line(0,1){10}}
\put(15,25){\line(1,0){10}}
\put(20,30){\line(1,0){25}}
\put(15,25){\line(1,3){5}}
\put(10,10){\line(1,3){5}}
\put(3,12){{$a_1$}}
\put(46,29){{$a_2$}}
\put(19,46){{$a_3$}}
\put(16,12){{$\Delta_{\epsilon_2}$}}
\put(30,42){{$\Delta_{\epsilon_1}$}}
\put(8,20){{$s$}}
\put(13,25){{$t$}}
\put(17,39){{$t$}}
\put(39,27){{$s$}}
\put(15,5){{Figure 2}}
\put(65,10){\line(1,1){20}}
\put(60,20){\line(1,0){15}}
\put(60,20){\line(4,1){40}}
\put(75,20){\line(1,2){10}}
\put(85,30){\line(0,1){15}}
\put(80,25){\line(0,1){5}}
\put(100,30){\line(0,1){10}}
\put(85,30){\line(1,0){20}}
\put(60,20){\line(2,1){40}}
\put(85,40){\line(1,0){15}}
\put(63,8){{$a_1$}}
\put(106,29){{$a_2$}}
\put(84,46){{$a_3$}}
\put(66,27){{$\Delta_{\epsilon_3}$}}
\put(95,42){{$\Delta_{\epsilon_1}$}}
\put(74,17){{$t$}}
\put(80,22){{$s$}}
\put(99,27){{$s$}}
\put(82,39){{$t$}}
\put(75,5){{Figure 3}}
\end{picture}

If we take $\alpha_2=sv_1+tv_2 \in \mathcal{C}_2$, then we have $t > s > 0$. Moreover, we have 
$\Theta_{cpt}(\alpha_2)= \{\epsilon_1, \epsilon_3 \}$, where
\begin{equation*}
\begin{split}
\epsilon_1 (i) & =1 \hspace{2mm} \text{for $i=1,2,3,4,5$} 
\end{split}, \hspace{5mm}
\begin{split}
\epsilon_3 (i) & =\left\{
\begin{array}{ll}
1 & \hspace{2mm} \text{for $i=1,3,5$} \\
-1 & \hspace{2mm} \text{for $i=2,4$.} \\
\end{array} \right.
\end{split}
\end{equation*}
The associated polyhedral complex $\mathcal{C}(X(\alpha_2,0))$ consists of 
all faces of $\Delta_{\epsilon_1}$ and $\Delta_{\epsilon_3}$ as in Figure 3. 

Since $\mathcal{C}(X(\alpha_1,0))$ and $\mathcal{C}(X(\alpha_2,0))$ have different 
combinatorial structures, the topology of the core changes. 
However, by Proposition \ref{indep}, $X(\alpha_1,0)$ and $X(\alpha_2,0)$ are diffeomorphic. 
Moreover the core is a deformation retract of the toric hyperk\"ahler manifold. 
Therefore we see that the cohomology ring of the core does not change, even if its topology 
changes. This is an interesting phenomenon. 

The variation of $\alpha$ corresponds to parallel translations of the hyperplanes $F_1, \dots F_N$. 
On the other hand, according to Theorem \ref{ring}, the cohomology ring of the core depends only on 
the subsets $\{ J \subset \{ 1, \dots, N \} \vert \bigcap_{i \in J} F_i = \emptyset \}$.
These subsets are preserved under parallel translations of the hyperplanes, 
even if the combinatorial structure of the associated polyhedral complex changes. 
This fact corresponds to the phenomenon mentioned above.}
\end{example}

\section{Variation of hyperk\"ahler structures}\label{variation}

In \cite{Ko3}, the author described how a toric hyperk\"ahler variety $X(\alpha,\beta)$ changes 
according to the variation of the parameter $(\alpha, \beta)$. 
First we can explain the geometric meaning of the parameter as follows. 
\begin{proposition}\label{period}
Suppose that the subtorus $K$ satisfies the condition that $X_i \not\in k$ for $i=1, \dots, N$. 
Then the following holds:
\newline
$(1)$ The hyperk\"ahler Kirwan map 
$\kappa \colon H_K^*(\mathbb{H}^N; \mathbb{R}) \to H^*(X(\alpha,\beta); \mathbb{R})$ 
induces an isomorphism $k^* \cong H^2(X(\alpha,\beta); \mathbb{R})$. 
\newline
$(2)$ The parameter $(\alpha,\beta)$ corresponds to the period of the 
hyperk\"ahler structure of $X(\alpha,\beta)$. That is , the following holds:
\begin{align*}
\kappa(\alpha)&=[\omega_1] \in 
H^2(X(\alpha,\beta);\mathbb{R}), \\
\kappa(\beta)&=
[\omega_{\mathbb{C}}] \in H^2(X(\alpha,\beta);\mathbb{C}).
\end{align*}
\end{proposition}
See \cite{Ko3} for the proof. $(1)$ is just a corollary of Theorem \ref{ring}.

Fix $\beta \in k^*_{\mathbb{C}}$. Let us describe how the toric hyperk\"ahler variety $X(\alpha,\beta)$ changes 
according to the variation of $\alpha$.  
The set $(k^* \times \{ \beta \}) \cap (k^* \times k^*_{\mathbb{C}})_{reg}$ is homeomorphic to $k^* \setminus \bigcup_{s \in S_{\beta}} W_s$, where $S_{\beta} = \{ s \in \{ 1, \dots, l \} ~|~ \beta \in W_{s} \otimes \mathbb{C} \}$.
A connected component of the set $k^* \setminus \bigcup_{s \in S_{\beta}} W_s$ is called a {\it chamber}. 
Note that the chamber structure depends on $\beta \in k^*_{\mathbb{C}}$. See Example \ref{ex2} for a 
concrete example. 

\begin{proposition}\label{isom}
Let $k^* \setminus \bigcup_{s \in S_{\beta}} W_s$ be the chamber structure for a 
fixed $\beta \in k^*_{\mathbb{C}}$. Suppose that $\alpha$ and $\alpha^{\prime}$ belong to the same chamber. 
Then $(X(\alpha^{\prime},\beta), I_1)$ is canonically biholomorphic to 
$(X(\alpha,\beta), I_1)$. 
\end{proposition}
\textit{Proof.} By Lemma \ref{criterion}, we have 
$\mu_{\mathbb{C}}^{-1}(\beta)^{\alpha -ss}=\mu_{\mathbb{C}}^{-1}(\beta)^{\alpha^{\prime} -ss}$.
So we finish the proof.  $\hfill \Box$
\vspace{3mm} \newline

Next we discuss the relation between two complex structures corresponding to different chambers. 
Namely, we study the so called wall-crossing phenomena. 

To do this, let us recall a special bimeromorphic map between complex symplectic manifolds, 
which was discovered by Mukai \cite{M}. Let $X_+$ be a complex symplectic manifold of complex dimension 
$2n$, containing a complex submanifold $V_+$ biholomorphic to $\mathbb{C}P^n$. 
Since $V_+$ is a complex Lagrangian submanifold, its normal bundle is biholomorphic to 
the cotangent bundle $T^*V_+$. Let $\pi \colon \widehat{X_+} \to X_+$ be the blowing-up along $V_+$. Then the 
exceptional divisor $\pi^{-1}(V_+)$ is biholomorphic to the projective cotangent bundle $\mathbb{P}(T^*V_+)$. 
This can be viewed as the variety of pairs $(p,h)$ such that $p \in V_+, h \in V_-$ with $p \subset h$, 
where $V_-$ is the space of hyperplanes in $V_+$. Therefore we can blow-down $\mathbb{P}(T^*V_+)$ to $V_-$ and 
get a new complex symplectic manifold $X_-$ and the bimeromorphic map $f \colon X_+ 
\dashrightarrow X_-$, which is called the \textit{Mukai flop} along $V_+$. 
This construction can be generalized directly to the case where $V_+$ is a $\mathbb{C}P^r$-bundle. 

\begin{example}
\rm{This is a continuation of Example \ref{ex1}. Suppose that $\beta=0$. 
Let $\pi_{\pm} \colon X(\alpha_{\pm},0) \to X(0,0)$ be the natural morphism. 
Then $V_0= \{ [0,0] \}$ is a set of singular points in $X(0,0)$. 
Note that $(0,0) \in \mathbb{H}^{n+1}$ is the unique fixed point for the action of $K_{\mathbb{C}}$. 
If we set $V_{\pm}=\pi_{\pm}^{-1}(V_0)$, then, by (\ref{plus}) and (\ref{minus}),  we have 
\begin{align*}
& V_+ = \{ (z,0) \in \mathbb{H}^{n+1} ~|~ z \ne 0 \} / K_{\mathbb{C}} \subset X(\alpha_+,0), \\
& V_- = \{ (0,w) \in \mathbb{H}^{n+1} ~|~ w \ne 0 \} / K_{\mathbb{C}} \subset X(\alpha_-,0),
\end{align*}
and $V_+ \cong V_- \cong \mathbb{C}P^n$. Therefore, by (\ref{plus}) and (\ref{minus}) again, we have 
\begin{equation*}
X(\alpha_+,0) \setminus V_+=X(\alpha_-,0) \setminus V_- 
=(\mu_{\mathbb{C}}^{-1}(0)^{\alpha_+ -ss} \cap \mu_{\mathbb{C}}^{-1}(0)^{\alpha_- -ss})/K_{\mathbb{C}}.
\end{equation*} 
Moreover, by Lemma \ref{criterion}, we see that, if $(z,w) 
\in \mu_{\mathbb{C}}^{-1}(0)^{\alpha_+ -ss} \cap \mu_{\mathbb{C}}^{-1}(0)^{\alpha_- -ss}$, 
then the orbit $(z,w)K_{\mathbb{C}}$ is closed in $\mu_{\mathbb{C}}^{-1}(0)$. Therefore 
$\pi_{\pm}|_{X(\alpha_{\pm},0) \setminus V_{\pm}}\colon X(\alpha_{\pm},0) \setminus V_{\pm} 
\to X(0,0) \setminus V_{0}$ is a biholomorphic map. 
Thus we have a bimeromorphic map $f \colon X(\alpha_+, 0) \dashrightarrow X(\alpha_-, 0)$ 
such that $f \vert_{X(\alpha_+,0) \setminus V_+} \colon X(\alpha_+,0) \setminus V_+ \to 
X(\alpha_-,0) \setminus V_-$ is a biholomorphic map. 
This is a typical example of wall-crossing phenomena and Mukai flops.}
\end{example}

General situations are described as follows. See \cite{Ko3} for more detailed discussions. 
\begin{theorem}\label{crossing} 
Suppose that the subtorus $K$ satisfies the condition \textnormal{(\ref{cmfd})}.
Let $k^* \setminus \bigcup_{s \in S_{\beta}} W_s$ be the chamber structure for a 
fixed $\beta \in k^*_{\mathbb{C}}$. Suppose that the chamber $\mathcal{C}_-$ is next to the chamber 
$\mathcal{C}_+$ across the wall $W_{s_0}$. 
Let $H^{(1)}_{s_0}$ be the 1-dimensional isotropy subgroup corresponding to the wall $W_{s_0}$. 
Set $J_{s_0}= \{ i \in \{ 1, \dots.N \} ~|~ \langle u_i, Y_{s_0} \rangle \ne 0 \}$, where 
$Y_{s_0}$ is a non-zero element in $Lie H^{(1)}_{s_0}$. 
Fix $\alpha_+ \in \mathcal{C}_+$ and $\alpha_- \in \mathcal{C}_-$. 
Fix $\alpha \in \overline{\mathcal{C}_+} \cap W_{s_0}$ such that $\alpha \not\in W_s$ 
for any $s \in S_{\beta}\setminus \{ s_0 \}$. Then the following holds: 
\newline
$(1)$ If we set $V_0= \{ [z,w] \in X(\alpha, \beta) ~|~ \text{ $(z,w)\zeta =(z,w)$ for $ \zeta \in 
H^{(1)}_{s_0}$} \}$, then $V_0$ is a toric hyperk\"ahler manifold. 
\newline
$(2)$ $\mu_{\mathbb{C}}^{-1}(\beta)^{\alpha_{\pm}-ss} \subset \mu_{\mathbb{C}}^{-1}(\beta)^{\alpha-ss}$ 
holds respectively. 
So we have the natural morphisms $\pi_{\pm} \colon (X(\alpha_{\pm},\beta),I_1) \to (X(\alpha,\beta),I_1)$.
\newline 
$(3)$ If we set $V_{\pm}= \pi_{\pm}^{-1}(V_0)$, then $\pi_{\pm}|_{V_{\pm}} \colon V_{\pm} \to V_0$ 
is a fiber bundle whose fiber is biholomorphic to $\mathbb{C}P^{\# J_{s_0}-1}$. 
Moreover, the codimension of $V_{\pm}$ in $X(\alpha_{\pm},\beta)$ is $\# J_{s_0}-1$, where $\# J_{s_0}$ is 
the number of elements in $J_{s_0}$.
\newline
$(4)$ The natural morphism $\pi_{\pm}|_{X(\alpha_{\pm},\beta) \setminus V_{\pm}} \colon X(\alpha_{\pm},\beta) 
\setminus V_{\pm} \to X(\alpha,\beta) \setminus V_0$ is a biholomorphic map. 
\end{theorem}
\textit{Proof.} $(1)$ Note that $W_{s_0}$ can be identified with the dual space of the Lie algebra $k_{s_0}$ of 
the quotient torus $K_{s_0}=K / H^{(1)}_{s_0}$. So $(\alpha,\beta)$ can be considered as an element 
of $k_{s_0}^* \times k_{s_0 \mathbb{C}}^*$. The assumption of $\alpha$ implies $(\alpha,\beta) \in 
(k_{s_0}^* \times k_{s_0 \mathbb{C}}^*)_{reg}$.
Then $V_0$ is a hyperk\"ahler quotient of $\mathbb{H}^{\{ 1, \dots, N \} \setminus J_{s_0}}= \{ (z,w) \in 
\mathbb{H}^N~|~ \text{$z_i=w_i=0$ if $i \in J_{s_0}$}\}$ by $K_{s_0}$ at $(\alpha,\beta) \in 
(k_{s_0}^* \times k_{s_0 \mathbb{C}}^*)_{reg}$. 
\newline
$(2)$ By Lemma \ref{criterion}, it is obvious. 
\newline
$(3)$ We can choose $Y_{s_0} \in Lie H^{(1)}_{s_0}$ so that $\langle \alpha_+, Y_{s_0} \rangle >0$. 
We set $J_{s_0}^{+} = \{ i ~|~ \langle u_i, Y_{s_0} \rangle >0 \}$ and 
$J_{s_0}^{-} = \{ i ~|~ \langle u_i, Y_{s_0} \rangle <0 \}$ respectively. 
By Lemma \ref{criterion} we can show that
\begin{equation}\label{ai}
\mu_{\mathbb{C}}^{-1}(\beta)^{\alpha_+ -ss}= \{ (z,w) \in \mu_{\mathbb{C}}^{-1}(\beta)^{\alpha -ss}
~|~ \text{(\ref{aa}) is satisfied} \},
\end{equation}
where
\begin{equation}\label{aa}
\text{there exists $i \in J_{s_0}$ such that $z_i \ne 0$ if $i \in J_{s_0}^+$ or $w_i \ne 0$ if 
$i \in J_{s_0}^-$.}
\end{equation}
It is also easy to see that, if $(z,w) \in \mu^{-1}(\alpha_+, \beta)$, then $[z,w] \in V_+$ is equivalent to 
\begin{equation}\label{b}
\text{$w_i=0$ for $i \in J_{s_0}^+$ and $z_i=0$ for $i \in J_{s_0}^-$.}
\end{equation}
Thus the fiber of $\pi_+|_{V_+} \colon V_+ \to V_0$ is biholomorphic to $(\mathbb{C}^{J_{s_0}} \setminus 
\{ 0 \})/H^{(1)}_{s_0\mathbb{C}}$. By the assumption that $X(\alpha_+,\beta)$ is non-singular, 
we see that the fiber is biholomorphic to $\mathbb{C}P^{\# J_{s_0}-1}$. 
Obviously the codimension of $V_+$ in $X(\alpha_+,\beta)$ is $\# J_{s_0}-1$.  
\newline
$(4)$ By Lemma \ref{criterion} we also have
\begin{equation}\label{c}
\mu_{\mathbb{C}}^{-1}(\beta)^{\alpha_- -ss}= \{ (z,w) \in \mu_{\mathbb{C}}^{-1}(\beta)^{\alpha -ss}
~|~ \text{(\ref{cc}) is satisfied} \},
\end{equation}
where
\begin{equation}\label{cc}
\text{there exists $i \in J_{s_0}$ such that $w_i \ne 0$ if $i \in J_{s_0}^+$ or $z_i \ne 0$ if 
$i \in J_{s_0}^-$.} 
\end{equation}
So (\ref{b}) implies
\begin{equation*}
X(\alpha_{+},\beta) \setminus V_+=
(\mu_{\mathbb{C}}^{-1}(\beta)^{\alpha_+ -ss} \cap \mu_{\mathbb{C}}^{-1}(\beta)^{\alpha_- -ss})/K_{\mathbb{C}}.
\end{equation*}
Similarly we have
\begin{equation*}
X(\alpha_{-},\beta) \setminus V_-
=(\mu_{\mathbb{C}}^{-1}(\beta)^{\alpha_- -ss} \cap \mu_{\mathbb{C}}^{-1}(\beta)^{\alpha_+ -ss})/K_{\mathbb{C}}.
\end{equation*}
On the other hand, by Lemma \ref{criterion}, if $(z,w) \in \mu_{\mathbb{C}}^{-1}(\beta)^{\alpha_+ -ss} 
\cap \mu_{\mathbb{C}}^{-1}(\beta)^{\alpha_- -ss}$, then the orbit $(z,w)K_{\mathbb{C}}$ is closed in 
$\mu_{\mathbb{C}}^{-1}(\beta)^{\alpha -ss}$. Thus we finish the proof. $\hfill \Box$
\vspace{3mm} \newline 
\begin{remark} \rm{In the above theorem we assumed the condition (\ref{cmfd}), that is, $X(\alpha_{\pm},\beta)$ 
is non-singular. Even if we drop this assumption, almost the same results hold. The exceptions 
are that $V_0$ is an orbifold and that the fiber of $\pi_{\pm}|_{V_{\pm}} \colon V_{\pm} \to V_0$ is a 
weighted projective space. The above proof works in this case. }
\end{remark}
\begin{proposition} Under the same assumption as Theorem \ref{crossing} we have the following.
\newline
$(1)$ If $\# J_{s_0} \ge 3$, then $(X(\alpha_-,\beta),I_1)$ is related to $(X(\alpha_+,\beta),I_1)$ by 
a Mukai flop. 
\newline
$(2)$ If $\# J_{s_0} =2$, there exists a biholomorphic map 
$\phi \colon \!(X(\alpha_+,\beta),I_1)\! \to \!(X(\alpha_-,\beta),I_1)$ satisfying $\pi_+ = \pi_- \circ \phi$.
\end{proposition} 
\textit{Proof.} $(1)$ By Theorem \ref{crossing} $(3)$ and $(4)$, we have a bimeromorphic map 
$$f \colon X(\alpha_{+},\beta) \dashrightarrow X(\alpha_{-},\beta)$$
such that $f |_{X(\alpha_{+},\beta) \setminus V_{+}} \colon X(\alpha_{+},\beta) \setminus 
V_{+} \to X(\alpha_{-},\beta) \setminus V_{-}$ is a biholomorphic map. By the proof of Theorem \ref{crossing} $(3)$, 
it is a fiber wise Mukai flop along $V_+$, which we also call a Mukai flop.  
\newline
$(2)$ We may assume $Y_{s_0}=X_1 \pm X_2$. In the following we assume $Y_{s_0}=X_1 + X_2$ and $\langle \alpha_+, Y_{s_0} \rangle >0$. 
(In the case $Y_{s_0}=X_1 - X_2$, a similar argument works.)

Define a map $\widetilde{\phi} \colon \mu_{\mathbb{C}}^{-1}(\beta)^{\alpha_+ -ss} \to 
\mu_{\mathbb{C}}^{-1}(\beta)^{\alpha_- -ss}$ by 
$$
\widetilde{\phi}((z,w))=((w_2,-w_1,z_3, \dots,z_N),(-z_2,z_1,w_3, \dots, w_N)).
$$
Suppose $(z,w) \in \mu_{\mathbb{C}}^{-1}(\beta)^{\alpha_+ -ss}$. 
Since $0= \langle \beta, Y_{s_0} \rangle=2 \pi \sqrt{-1}(z_1w_1+z_2w_2)$, we have 
$\mu_{\mathbb{C}}(\widetilde{\phi}(z,w))=\mu_{\mathbb{C}}(z,w)=\beta$. 
By (\ref{ai}) we have $(z_1,z_2) \ne (0,0)$. Moreover, it is easy to see the following;
\begin{equation}
\begin{split}\label{sss}
&\text{if $(z,w) \in \mu_{\mathbb{C}}^{-1}(\beta)^{\alpha_+ -ss}$ and $(w_1,w_2) \ne (0,0)$, }\\
&\text{\hspace*{10mm}then there exists $t \in \mathbb{C}$ such that $\widetilde{\phi}((z,w))=(z,w)\mathrm{Exp}tY_{s_0}$. }
\end{split}
\end{equation}
Moreover, if $(z,w) \in \mu_{\mathbb{C}}^{-1}(\beta)^{\alpha_+ -ss}$ and $(w_1,w_2) = (0,0)$, we have
$$\alpha \in \sum_{i=3}^N \mathbb{R}_{\ge 0}|z_i|^2\iota^{*}u_i
+\sum_{i=3}^N \mathbb{R}_{\ge 0}|w_i|^2 (-\iota^{*}u_i).
$$
Thus we see that $\widetilde{\phi}((z,w)) \in \mu_{\mathbb{C}}^{-1}(\beta)^{\alpha -ss}$.
Moreover, by (\ref{c}), we have $\widetilde{\phi}((z,w)) \in \mu_{\mathbb{C}}^{-1}(\beta)^{\alpha_- -ss}$. 
Therefore the map $\widetilde{\phi}$ is well-defined.

Note that the map $\widetilde{\phi}$ is \textit{not} $K_{\mathbb{C}}$-equivariant. 
However, it is easy to see that $\widetilde{\phi}((z,w)\mathrm{Exp}tY_{s_0})=(z,w)\mathrm{Exp}(-tY_{s_0})$. 
Moreover, if $Z \in k$ satisfies $\langle \iota^* u_1, Z \rangle=1$ and 
$\langle \iota^* u_2, Z \rangle=-1$, then we have $\widetilde{\phi}((z,w)\mathrm{Exp}tZ)=(z,w)\mathrm{Exp}tZ$. 
If $Z \in k$ satisfies $\langle \iota^* u_1, Z \rangle=\langle \iota^* u_2, Z \rangle=0$, 
then we have $\widetilde{\phi}((z,w)\mathrm{Exp}tZ)=(z,w)\mathrm{Exp}tZ$. Thus $\widetilde{\phi}$ 
induces a biholomorphic map $\phi \colon (X(\alpha_+,\beta),I_1) \to (X(\alpha_-,\beta),I_1)$.

Finally, by (\ref{sss}) again, we have $\pi_+ = \pi_- \circ \phi$. $\hfill \Box$
\vspace{3mm} \newline 
Variation of GIT quotients was studied in \cite{DH} and \cite{T1}. 
Our cases are more restrictive, because we are treating with hyperk\"ahler quotients. 
So the bimeromorphic map between two quotients is a special one, that is, a Mukai flop. 
Similar phenomena were observed in \cite{T2}.

In the above proof we saw that there exists a stratification $V_0 \subset X(\alpha,\beta)$,  
which comes from isotropy subgroups. More generally, if $(\alpha,\beta) \in k^* \times k_{\mathbb{C}}^* \setminus 
(k^* \times k_{\mathbb{C}}^*)_{reg}$, then a point $(z,w) \in \mu_{\mathbb{C}}^{-1}(\beta)$ may have 
an isotropy subgroup $K_{(z,w)} \subset K$ of positive dimension. So the toric hyperk\"ahler variety 
$X(\alpha,\beta)$ has a stratification, which is indexed by the set of isotropy subgroups.  
Proudfoot and Webser used this stratification to compute the intersection cohomology of 
singular toric hyperk\"ahler varieties \cite{PW}. 

\vspace{3mm}
\noindent
\footnotesize{Graduate School of Mathematical Sciences, The University of Tokyo, \\
3-8-1, Komaba, Meguro-ku, Tokyo, 153-8914, Japan \\ 
E-mail address: konno@ms.u-tokyo.ac.jp}
\end{document}